\documentclass[dvipsnames]{amsart}
\usepackage{pict2e}
\usepackage{amsmath,amssymb,amsthm}
\usepackage{mathrsfs}
\usepackage{fullpage}
\usepackage{todonotes}

\allowdisplaybreaks
\usepackage{soul}

\usepackage{color}
\usepackage{hyperref}
\hypersetup{
	colorlinks=true,
	linkcolor=blue,
	citecolor=blue,
	urlcolor=blue
}

\usepackage{xcolor}

\usepackage{tikz}

\newtheorem{theorem}{Theorem}[section]
\newtheorem{lemma}[theorem]{Lemma}	
\newtheorem{proposition}[theorem]{Proposition}
\newtheorem{corollary}[theorem]{Corollary}
\theoremstyle{definition}
\newtheorem{definition}[theorem]{Definition} 
\newtheorem{remark}[theorem]{Remark}	
\newtheorem{example}[theorem]{Example}

 \usepackage{relsize}
\theoremstyle{definition} 

\usetikzlibrary{
  arrows.meta,
  calc,
  decorations.markings,
  decorations.pathmorphing,
  intersections,
  patterns,
  positioning,
  angles,
  quotes
}
\usepackage{todonotes}
\usepackage{booktabs}
\usepackage{tikz-cd}
\usetikzlibrary{positioning, angles}
\usetikzlibrary{arrows}
\usetikzlibrary{decorations.pathmorphing, decorations.pathreplacing}
\usetikzlibrary{decorations.markings}
\usepackage{mathtools}
\usepackage{tikz}
\newtheorem{introtheorem}{Theorem}

\usepackage{tikz-3dplot}
\usetikzlibrary{3d, angles, arrows.meta, calc, decorations.markings, 
                decorations.pathmorphing, decorations.pathreplacing, 
                positioning, shapes,shadings}
\usepackage{pgfplots}
\usepackage{pgfmath}

\usepackage[nameinlink,noabbrev]{cleveref}
\usepackage{verbatim}
\theoremstyle{plain}

\theoremstyle{remark}

\newtheorem*{theorem*}{Theorem}
\newtheorem*{remark*}{Remark}

\crefname{theorem}{theorem}{theorems}
\Crefname{theorem}{Theorem}{Theorems}
\crefname{lemma}{lemma}{lemmas}
\Crefname{lemma}{Lemma}{Lemmas}
\crefname{proposition}{proposition}{propositions}
\Crefname{proposition}{Proposition}{Propositions}
\crefname{corollary}{corollary}{corollaries}
\Crefname{corollary}{Corollary}{Corollaries}
\crefname{conjecture}{conjecture}{conjectures}
\Crefname{conjecture}{Conjecture}{Conjectures}
\Crefname{construction}{Construction}{Construction}
\crefname{definition}{definition}{definitions}
\Crefname{definition}{Definition}{Definitions}
\crefname{example}{example}{examples}
\Crefname{example}{Example}{Examples}
\crefname{exercise}{exercise}{exercises}
\Crefname{exercise}{Exercise}{Exercises}
\crefname{assumption}{assumption}{assumptions}
\Crefname{assumption}{Assumption}{Assumptions}
\crefname{question}{question}{questions}
\Crefname{question}{Question}{Questions}
\crefname{remark}{remark}{remarks}
\Crefname{remark}{Remark}{Remarks}
\crefname{notation}{notation}{notations}
\Crefname{notation}{Notation}{Notations}
\crefname{equation}{equation}{equations}
\Crefname{equation}{Equation}{Equations}
\crefname{figure}{figure}{figures}
\Crefname{figure}{Figure}{Figures}
\crefname{table}{table}{tables}
\Crefname{table}{Table}{Tables}

\author[C. Sopio]{Cristian Sopio}  
\address{Cristian Sopio\\Dipartimento di Matematica e Informatica, Università degli Studi di Ferrara, Via Machiavelli 35, 44121 Ferrara, Italy}  \email{\url{cristian.sopio@unife.it}}
\title{An intrinsic criterion for Pansu differentials to be homomorphisms of nilpotent Lie groups}
\author[F. Tripaldi]{Francesca Tripaldi} 
\address{Francesca Tripaldi\\Department of Pure Mathematics, University of Leeds,  Woodhouse, LS2 9JT Leeds, UK}
\email{\url{f.tripaldi@leeds.ac.uk}}
\subjclass{53C17, 22E25, 17B70}

\keywords{Pansu differentiability, Carnot groups, nilpotent Lie groups, adapted gradings, osculating groups, Lipschitz maps}
\begin{document}

\maketitle
\begin{abstract}
By Pansu's theorem, Lipschitz maps between Carnot groups are
differentiable almost everywhere, and their differentials are group
homomorphisms. When the domain is instead a nilpotent Lie group with a
bracket-generating left-invariant distribution, the Pansu differential
is a homomorphism of an osculating Carnot group, which depends on an
auxiliary choice of grading, and it need not be a homomorphism of the
original group. We give an intrinsic criterion for when it is one, for
a suitable choice of grading: the differential must vanish on an ideal
that depends only on the group and its distribution, and comes from a
universal quotient onto a Carnot group. The criterion holds
automatically when the step of the target is small enough, and an
explicit example shows that this threshold is sharp. Along the way, we
determine exactly how the obstruction depends on the choice of grading.
\end{abstract}
\section{Introduction}

Consider a nilpotent Lie group $G$ with a left-invariant
bracket-generating distribution, and a Carnot group $H$. Lipschitz maps
from $G$ to $H$ are Pansu differentiable almost everywhere \cite{KV08}.
When $G$ is itself a Carnot group whose first layer is the given
distribution, by Pansu's theorem \cite{Pan89} we know that the differentials
are group homomorphisms $G\to H$. In general, the Pansu differential is
a homomorphism from the tangent cone of $G$, which is a Carnot group
that need not be isomorphic to $G$ (Mitchell's theorem \cite{Bel96}, see also
Theorem~\ref{thm:mitchell}). Moreover, it can be written as a map on
$G$ only after choosing a grading adapted to the distribution. We give
an intrinsic criterion for when, for a suitable choice of grading, the
Pansu differential is a group homomorphism between the two given groups $G\to H$, and we describe
exactly how the obstruction depends on the choice of grading.

Let $G$ be a connected, simply connected nilpotent Lie group with Lie
algebra $\mathfrak g=(E,[\cdot,\cdot]_{\mathfrak g})$, where $E$ is a
finite-dimensional real vector space. A bracket-generating subspace
$\Delta\subseteq E$ defines an increasing filtration
$(\mathcal F_iE)_{i\geq 0}$: we set $\mathcal F_0E=\{0\}$ and, for
$i\geq 1$, let $\mathcal F_iE$ be the linear span of bracket words of
length at most $i$ in $\Delta$. Thus $\mathcal F_1E=\Delta$ and
$[\mathcal F_iE,\mathcal F_jE]_{\mathfrak g}\subseteq\mathcal F_{i+j}E$.
The \textit{distribution step} is the smallest integer $s$ such that
$\mathcal F_sE=E$, and in general $s\le S$, where $S$ is the nilpotency
step of $\mathfrak g$. The associated graded Lie algebra
\begin{equation}\label{eq:mathfrak-v-defined}
 \mathfrak v:=\operatorname{gr}_{\Delta}\mathfrak g
 =\bigoplus_{j=1}^{s}\mathfrak v_j\ ,
 \quad
 \mathfrak v_j=\mathcal F_jE/\mathcal F_{j-1}E\,,
\end{equation}
with the bracket induced by $[\cdot,\cdot]_{\mathfrak g}$, is a Carnot
algebra with first layer $\mathcal{F}_1E=\Delta$. We call
$\mathfrak v_j$ the layer of \textit{weight $j$}, and a nonzero element
of $\mathfrak v_j$ is called \emph{homogeneous of weight $j$}. We use
the same terminology for other graded vector spaces.

Let $H$ be a Carnot group with stratified Lie algebra
$\mathfrak h=\bigoplus_{j=1}^{m}\mathfrak h_j$ of step $m$. We write
$\mathfrak h_{\leq i}=\bigoplus_{j=1}^{i}\mathfrak h_j$, with
$\mathfrak h_j=\{0\}$ for $j>m$, and use the same convention for other
graded algebras. A Lie group homomorphism $\Phi\colon G\to H$ is
\emph{horizontal} if its differential at the identity
$d\Phi_e\colon\mathfrak g\to\mathfrak h$ satisfies
$d\Phi_e(\Delta)\subseteq\mathfrak h_1$. Consequently,
$d\Phi_e(\mathcal{F}_iE)\subseteq\mathfrak h_{\le i}$, and so
$d\Phi_e$ induces a \textit{graded} Lie algebra homomorphism
$\operatorname{gr}d\Phi_e\colon\mathfrak v\to\mathfrak h$ mapping each
$\mathfrak v_j$ into $\mathfrak h_j$. The question we address is:
\begin{center}
Which graded Lie algebra homomorphisms
$D\colon\mathfrak v\to\mathfrak h$ have the form
$D=\operatorname{gr}d\Phi_e$ \\for a horizontal Lie group homomorphism
$\Phi\colon G\to H$?
\end{center}
Since $G$ is simply connected, this is equivalent to asking whether
$D=\operatorname{gr}L$ for some Lie algebra homomorphism
$L\colon\mathfrak g\to\mathfrak h$ satisfying
$L(\Delta)\subseteq\mathfrak h_1$. If such an $L$, and hence $\Phi$,
exists, it is unique, since $L|_{\Delta}=D|_{\mathfrak v_1}$ and
$\Delta$ bracket-generates $\mathfrak g$
(Lemma~\ref{lem:horizontal-gr}).

\smallskip
\textbf{The universal horizontal quotient.}
Proposition~\ref{prop:horizontal-quotient} constructs a Carnot algebra
$\mathfrak G$ and a surjective Lie algebra homomorphism
$\rho\colon\mathfrak g\to\mathfrak{G}$ with $\rho(\Delta)=\mathfrak G_1$,
with the following universal property: every Lie algebra homomorphism
$L\colon\mathfrak g\to\mathfrak h$ into a Carnot algebra with
$L(\Delta)\subseteq\mathfrak h_1$ factors uniquely as
$L=\overline{L}\circ\rho$, where
$\overline{L}\colon\mathfrak G\to\mathfrak h$ is graded:
\[
\begin{tikzcd}[column sep=large, row sep=large]
\mathfrak g
  \arrow[rr, two heads, "\rho"]
  \arrow[dr, "L"']
&&
\mathfrak G
  \arrow[dl, dashed, "{\exists!\,\overline L\ \text{graded}}"]
\\
& \mathfrak h &
\end{tikzcd}
\]
Thus $\ker\rho$ consists of the elements annihilated by every such $L$,
and $(\mathfrak G,\rho)$ depends only on $(\mathfrak g,\Delta)$, up to a
unique graded isomorphism intertwining the quotient maps
(Corollary~\ref{cor:horizontal-quotient}). In practice, $\mathfrak G$
is computed from a presentation of $\mathfrak g$ with generators in
$\Delta$, by replacing each relation with its homogeneous components.
Moreover, $\rho$ is strict, that is,
$\rho(\mathcal{F}_iE)=\mathfrak G_{\le i}$ for every $i$, so its
associated graded homomorphism
$\operatorname{gr}\rho\colon\mathfrak v\to\mathfrak G$ maps
$\mathfrak v_j$ onto $\mathfrak G_j$ for every $j$.

\smallskip
\textbf{Adapted gradings and defect ideals.} To compare $\mathfrak G$
with realisations of $\mathfrak v$ on $E$, we choose an \textit{adapted
grading} $W=(W_j)_{j=1}^s$, that is, a grading of the vector space $E$
with
\begin{equation*}
 E=\bigoplus_{j=1}^{s}W_j\ ,
 \quad
 \mathcal F_jE=\bigoplus_{i=1}^{j}W_i\ ,
 \quad W_1=\Delta\,,
\end{equation*}
called a tangent grading in \cite{LDNGT}. We do not require
$[W_i,W_j]_{\mathfrak g}\subseteq W_{i+j}$. The linear isomorphism
$\sigma_W\colon\mathfrak v\to E$ sends
$\xi\in\mathfrak v_j=\mathcal F_jE/\mathcal F_{j-1}E$ to its unique
representative in $W_j$. Transporting the bracket of $\mathfrak v$ via
$\sigma_W$ gives the \emph{osculating algebra}
$\mathfrak g_0^W=(E,[\cdot,\cdot]_{0,W})$, where
\[
 [u,v]_{0,W}
 :=\sigma_W\bigl([\sigma_W^{-1}(u),\sigma_W^{-1}(v)]_{\mathfrak v}\bigr)\ ,
 \quad u,v\in E\,.
\]
The \textit{bracket defect}
\[
 B_W(u,v):=[u,v]_{\mathfrak g}-[u,v]_{0,W}\ ,
 \quad u,v\in E\,,
\]
measures the disagreement between the original and osculating brackets,
and the \textit{defect ideal}
\begin{equation}\label{eq:defect-ideal}
 J_W:=\operatorname{Ideal}_{0,W}
 \bigl(\{\operatorname{pr}_{W_i}B_W(u,v):
 u,v\in E,\ 1\leq i\leq s\}\bigr)
\end{equation}
is the smallest homogeneous ideal of $\mathfrak g_0^W$ containing every
bracket defect. Here $\operatorname{pr}_{W_i}$ denotes the projection
onto $W_i$ along the other layers. To compare different gradings, we
transport $J_W$ to the fixed algebra $\mathfrak v$: the subspace
$j_W:=\sigma^{-1}_W(J_W)$ is a homogeneous ideal of $\mathfrak v$.

For a graded homomorphism $D\colon\mathfrak v\to\mathfrak h$, the map
$A_W:=D\circ\sigma_W^{-1}$ is a graded homomorphism from
$\mathfrak g_0^W$ to $\mathfrak h$, and
\[
 A_W([u,v]_{\mathfrak g})-[A_W(u),A_W(v)]_{\mathfrak h}
 =A_W(B_W(u,v))\,.
\]
Since $\ker A_W$ is a homogeneous ideal of $\mathfrak g_0^W$, this gives
the following criterion (Lemma~\ref{lem:kernel-criterion}):
\begin{equation}\label{eq:kernel-criterion-intro}
 A_W\text{ preserves }[\cdot,\cdot]_{\mathfrak g}
 \ \Longleftrightarrow\
 J_W\subseteq\ker A_W
 \ \Longleftrightarrow\
 j_W\subseteq\ker D\,.
\end{equation}
Moreover, $J_W$ is also an ideal of $\mathfrak g$, and
$\mathfrak q_W:=E/J_W$ is a common Carnot quotient of $\mathfrak g$ and
$\mathfrak g_0^W$ on which the two brackets agree
(Proposition~\ref{prop:common-quotient}).

\smallskip
\textbf{A tilted Heisenberg grading.}
Let $\mathfrak g$ be the three-dimensional Heisenberg algebra with basis
$X,Y,Z$ and $[X,Y]_{\mathfrak g}=Z$, and let
$\Delta=\operatorname{span}\{X,Y\}$. Then $\mathfrak v_1=\Delta$ and
$\mathfrak v_2=\operatorname{span}\{\overline Z\}$, where
$\overline Z:=Z+\Delta$. Both
\[
 W_1=\Delta\ ,\ W_2=\operatorname{span}\{Z\}
 \quad\text{and}\quad
 \widetilde{W}_1=\Delta\ ,\ \widetilde{W}_2=\operatorname{span}\{Z+X\}
\]
are adapted gradings, but only the first is a Lie algebra grading of
$\mathfrak g$. The maps $\sigma_W$ and $\sigma_{\widetilde W}$ agree on
$\mathfrak v_1$, while $\sigma_{W}(\overline Z)=Z$ and
$\sigma_{\widetilde{W}}(\overline Z)=Z+X$. Thus the associated graded
algebra is the same, but its realisation on $E$ depends on the adapted
grading.

For $W$, the brackets $[\cdot,\cdot]_\mathfrak g$ and
$[\cdot,\cdot]_{0,W}$ coincide, so $J_W=\{0\}$. For $\widetilde{W}$, in
the basis $X,Y,T$ with $T:=Z+X$, we have
\[
\begin{aligned}
 &[X,Y]_{\mathfrak g}=T-X\ ,
 &&[X,T]_{\mathfrak g}=0\ ,
 &&[Y,T]_{\mathfrak g}=X-T\ ,\\
 &[X,Y]_{0,\widetilde{W}}=T\ ,
 &&[X,T]_{0,\widetilde{W}}=0\ ,
 &&[Y,T]_{0,\widetilde{W}}=0\,.
\end{aligned}
\]
The defect $B_{\widetilde{W}}(X,Y)=-X$ puts $X$, and hence
$T=[X,Y]_{0,\widetilde{W}}$, in $J_{\widetilde{W}}$, so
$J_{\widetilde{W}}=\operatorname{span}\{X,T\}=\operatorname{span}\{X,Z\}$.
Writing $\overline{X}\in\mathfrak v_1$ for $X$, we get
\[
 j_W=\{0\}\quad\text{and}\quad
 j_{\widetilde{W}}=\operatorname{span}\{\overline{X},\overline{Z}\}\,.
\]
The same algebra $\mathfrak v$ thus carries different defect ideals. By
\eqref{eq:kernel-criterion-intro}, for an injective graded homomorphism
$D\colon\mathfrak v\to\mathfrak h$, the map $D\circ\sigma_W^{-1}$
preserves the original bracket, while $D\circ\sigma_{\widetilde{W}}^{-1}$
does not.

\smallskip
\textbf{Minimal gradings and the answer.} To relate the choice of
adapted grading to the universal quotient, consider the two linear maps
\[
 \operatorname{gr}\rho\colon\mathfrak v\to\mathfrak G
 \quad\text{and}\quad
 \rho\circ\sigma_W\colon\mathfrak v\to\mathfrak G\,.
\]
For $\xi\in\mathfrak v_j$, the vectors $\rho(\sigma_W(\xi))$ and
$\operatorname{gr}\rho(\xi)$ have the same component of weight $j$,
since $\sigma_W(\xi)$ represents $\xi$, while the components of lower
weight may depend on $W$. These are recorded by the \textit{correction
map}
\[
 t_W:=\rho\circ\sigma_W-\operatorname{gr}\rho\colon\mathfrak v\to\mathfrak G\ ,
 \quad t_W(\mathfrak v_j)\subseteq\mathfrak G_{\le j-1}\,.
\]
In the Heisenberg example, $\mathfrak G=\mathfrak g$ and
$\rho=\mathrm{id}_E$, as for any Carnot algebra with first layer
$\Delta$, and $t_{\widetilde{W}}(\overline{Z})=X$. Let $M_W$ be the
homogeneous ideal of $\mathfrak G$ generated by the homogeneous
components of the vectors $t_W(\xi)$, $\xi\in\mathfrak v$.
Proposition~\ref{prop:correction-ideal} gives
\begin{equation}\label{eq:comparison-intro}
 J_W=\rho^{-1}(M_W)\ ,
 \quad j_W=(\operatorname{gr}\rho)^{-1}(M_W)\ ,
 \quad \mathfrak q_W\cong\mathfrak G/M_W\,.
\end{equation}
We call $W$ \textit{minimal} if $t_W=0$, that is, if
$\rho\circ\sigma_W=\operatorname{gr}\rho$. Minimal gradings exist, and
by \eqref{eq:comparison-intro} they are exactly the adapted gradings
with the smallest possible defect ideals, namely $J_W=\ker\rho$ and
$j_W=\ker(\operatorname{gr}\rho)$. They need not be unique: they are
unique if and only if $\ker(\operatorname{gr}\rho)\subseteq\mathfrak v_s$
(Corollary~\ref{cor:least-ideals} and Remark~\ref{rem:nonunique}). We
can now answer our question.

\begin{introtheorem}[see Theorem~\ref{thm:intrinsic}]
\label{thm:A}
For a graded Lie algebra homomorphism
$D\colon\mathfrak v\to\mathfrak h$, the following are equivalent:
\begin{enumerate}
 \item[(i)] $D=\operatorname{gr}d\Phi_e$ for a horizontal Lie group
       homomorphism $\Phi\colon G\to H$;
 \item[(ii)] $D\circ\sigma_W^{-1}$ preserves $[\cdot,\cdot]_{\mathfrak g}$
       for some adapted grading $W$;
 \item[(iii)] $D\circ\sigma_W^{-1}$ preserves $[\cdot,\cdot]_{\mathfrak g}$
       for every minimal adapted grading $W$;
 \item[(iv)] $D(\ker\operatorname{gr}\rho)=\{0\}$.
\end{enumerate}
In that case, $\Phi$ is unique and $D\circ\sigma_W^{-1}=d\Phi_e$ for
every adapted grading $W$ as in \textup{(ii)}.
\end{introtheorem}

Thus a single minimal adapted grading, chosen using only
$(\mathfrak g,\Delta)$, realises $d\Phi_e$ for every such $D$ and every
Carnot target $H$, and the minimal gradings are the only adapted
gradings with this property. For a fixed $D$, Theorem~\ref{thm:intrinsic}
also describes all the adapted gradings $W$ for which
$D\circ\sigma_W^{-1}$ preserves $[\cdot,\cdot]_{\mathfrak g}$ (see the
discussion following it).

\smallskip
\textbf{Which defect ideals occur.} For a homogeneous ideal
$M\subseteq\mathfrak G$, let $[\mathfrak G,M]$ be the span of the
brackets of elements of $\mathfrak G$ with elements of $M$, and let
\[
 \mu_j(M):=\dim\bigl((M/[\mathfrak G,M])_j\bigr)\,.
\]
By Lemma~\ref{lem:minimal-generators}, $\mu_j(M)$ is the least number of
elements of weight $j$ in a set of homogeneous generators of $M$ as an
ideal. Let $\mathfrak v_{>j}:=\bigoplus_{i>j}\mathfrak v_i$.

\begin{introtheorem}[see Theorem~\ref{thm:attainability}]\label{thm:B}
A homogeneous ideal $M\subseteq\mathfrak G$ equals $M_W$ for some
adapted grading $W$ if and only if
$\mu_j(M)\leq\dim\mathfrak v_{>j}$ for every $j\geq1$.
\end{introtheorem}

The proof constructs an adapted grading for each such $M$. Together with
\eqref{eq:comparison-intro}, this determines all the defect ideals $J_W$
and $j_W$ of the pair $(\mathfrak g,\Delta)$. Since
$\mathfrak v_{>s}=\{0\}$, the ideal $M_W$ never needs a generator of the
top weight $s$. For instance, the centre of the Heisenberg algebra is
never a defect ideal (Example~\ref{ex:heisenberg-all}).
Proposition~\ref{prop:n51} lists all the defect ideals of a
five-dimensional nonstratifiable example.

\smallskip
\textbf{Pansu differentials and the target step.} Equip $G$ and $H$
with the left-invariant Carnot--Carath\'eodory distances determined by
inner products on $\Delta$ and $\mathfrak h_1$, and let
$f\colon\Omega\to H$ be Lipschitz on an open set $\Omega\subseteq G$.
For every adapted grading $W$, the map $f$ is Pansu differentiable with
respect to $W$ almost everywhere (Theorem~\ref{thm:existence}): at
almost every $x\in\Omega$, the blow-ups of $f$ at $x$, taken with the
dilations of $E$ that multiply $W_j$ by $t^j$, converge to a graded
homomorphism $D^W_Pf(x)\colon\mathfrak g_0^W\to\mathfrak h$, the
\emph{Pansu differential} of $f$ at $x$ with respect to $W$
(Definition~\ref{def:pansu-W}).

These differentials are representatives of one intrinsic object: there
is a set $\Omega_f\subseteq\Omega$ of full measure and, at every
$x\in\Omega_f$, a graded homomorphism
$D_Pf(x)\colon\mathfrak v\to\mathfrak h$, the \emph{canonical Pansu
differential}, such that $D^W_Pf(x)=D_Pf(x)\circ\sigma_W^{-1}$ for every
adapted grading $W$ (Lemma~\ref{lem:change-grading} and
Definition~\ref{def:canonical-differential}). The relevant invariant of
$(\mathfrak g,\Delta)$ is its \emph{visibility weight}, the lowest weight
of the homogeneous ideal $\ker(\operatorname{gr}\rho)$,
\[
 r_*:=\min\bigl\{j\geq1:
 \ker(\operatorname{gr}\rho)\cap\mathfrak v_j\neq\{0\}\bigr\}\,,
\]
with $r_*=\infty$ when $\operatorname{gr}\rho$ is injective
(Definition~\ref{def:rstar}). Theorem~\ref{thm:A}, applied to
$D=D_Pf(x)$, gives the following.

\begin{introtheorem}[see Theorem~\ref{thm:intrinsic-analytic} and
Corollary~\ref{cor:visibility}]\label{thm:C}
Let $f\colon\Omega\to H$ be Lipschitz and $x\in\Omega_f$. There is an
adapted grading $W$ for which $D^W_Pf(x)$ preserves
$[\cdot,\cdot]_{\mathfrak g}$, that is, $D^W_Pf(x)\colon G\to H$ is a
group homomorphism, if and only if
$D_Pf(x)(\ker\operatorname{gr}\rho)=\{0\}$. In that case, every minimal
adapted grading has this property, and $D^W_Pf(x)=d\Phi_e$ for a unique
horizontal homomorphism $\Phi\colon G\to H$. If $m<r_*$, this holds at
every $x\in\Omega_f$, for every Lipschitz map $f$.
\end{introtheorem}

The last assertion holds because $\ker(\operatorname{gr}\rho)$ is
homogeneous and lies in weights at least $r_*$, while $D_Pf(x)$ maps
$\mathfrak v_j$ into $\mathfrak h_j=\{0\}$ for $j>m$: the obstruction is
invisible to targets of step less than $r_*$. In particular, the
grading can then be fixed independently of $f$ and $x$.

For $m\geq r_*$, the criterion $D_Pf(x)(\ker\operatorname{gr}\rho)=\{0\}$
is a genuine condition, and the bound $m<r_*$ is sharp. For the
five-dimensional algebra $N_{5,1}$ with its natural two-dimensional
$\Delta$, the quotient $\mathfrak G$ is the Engel algebra and $r_*=3$,
and Proposition~\ref{prop:n51-engel} exhibits a polynomial contact map
into the Engel group, of step $3$, whose Pansu differential violates the
criterion at every point where it exists. The bound is not always
necessary, though. For the six-dimensional pair of
Example~\ref{ex:rigidity-algebra}, we have $r_*=2$, and for maps into
the first Heisenberg group, so that $m=r_*$, the criterion holds at
every point whenever the Pansu differential with respect to the minimal
grading exists everywhere and depends continuously on the point
(Proposition~\ref{prop:rigidity}). We do not know whether the same
conclusion holds almost everywhere for every Lipschitz map. By
contrast, in the Heisenberg example above $\operatorname{gr}\rho$ is
injective, so $r_*=\infty$, and the failure for $\widetilde W$ concerns
only that representative (Example~\ref{ex:heisenberg-revisited}).

Finally, Theorem~\ref{thm:A} holds verbatim when $\mathfrak h$ is only
graded, and the pointwise statements of Theorem~\ref{thm:C} hold for
graded targets at every point of Pansu differentiability. What is
missing for graded targets is the almost everywhere existence of Pansu
differentials when the domain is not a Carnot group
(Remark~\ref{rem:graded-targets}).

\smallskip
Adapted gradings appear as tangent gradings in
\cite{LDNGT}, where quantitative versions of Mitchell's theorem and of
Pansu's asymptotic theorem are proved. Our common quotient
$\mathfrak q_W$ is closely related to their Carnot quotient ideals
\cite[Definition~3.18 and Lemma~3.23]{LDNGT}: the defect ideal $J_W$ is
the smallest of them (Corollary~\ref{cor:least-compatible}). Osculating
groups go back to Folland and Stein \cite{FS74}, and for general filtered
manifolds, see \cite{vEY17}. Comparisons between a nilpotent Lie algebra
and its associated graded algebra also appear in \cite{Cor19}, where
gradings adapted to the lower central series yield quantitative
sublinear biLipschitz comparisons with the associated Carnot group, and in \cite{Tri20}, where positive gradings
are used to construct the Rumin complex. Our
focus is instead the dependence of the defect ideals on the adapted
grading, their classification, and the resulting intrinsic criterion.

The application to Lipschitz maps relies on almost everywhere Pansu
differentiability, due to Pansu \cite{Pan89} when $G$ is a Carnot group
(see also \cite[Theorem~2.1.6]{Monti2001}) and to Karmanova and
Vodopyanov \cite{KV08} in general. In the direction of the target,
Pansu's theorem has been extended to graded groups
\cite[Theorem~2.11]{Mag10} and to Banach homogeneous groups
\cite[Theorem~1.1]{MR14}. When the domain is a stratified group, locally
Lipschitz smooth maps are characterised by contact equations
\cite[Theorem~1.1]{Mag10}, and differential calculus in stratified
groups is studied in \cite{Mag13}. The Engel group, the target of our
sharpness example, is the simplest Carnot group of step $3$. The algebras
$N_{5,1}$ and $N_{6,2,10}$ in our examples are named as in the list of
low-dimensional nilpotent Lie algebras of \cite{LDT22}.

\smallskip

Section~\ref{sec:preliminaries} collects notation
and preliminaries on graded algebras, adapted gradings, defect ideals,
Carnot-Carath\'eodory distances and Pansu differentiability, and ends
with the change of grading and the canonical Pansu differential.
Section~\ref{sec:quotients} constructs the universal horizontal quotient
and the common quotient for a fixed adapted grading, defines the
visibility weight, and relates the two quotients through the correction
map. Section~\ref{sec:classification} describes all adapted gradings in
terms of a minimal one, proves Theorem~\ref{thm:B}, and determines all
defect ideals for the Heisenberg algebra and for $N_{5,1}$.
Section~\ref{sec:intrinsic} proves Theorems~\ref{thm:A}
and~\ref{thm:C}, and contains the sharpness and rigidity examples.
\subsection*{Acknowledgements}
The first author thanks the School of Mathematics at the University of Leeds for the hospitality during his research visit. The authors thank Andrea Tettamanti for bringing the paper \cite{LDNGT} to their attention.
\subsection*{Use of AI tools}
The authors used ChatGPT (OpenAI) to assist with the revision of the
manuscript and to help check the mathematical arguments.
The authors take responsibility for the final text.

\section{Preliminaries}
\label{sec:preliminaries}
 
We keep the notation of the introduction. Lie algebras are assumed
finite-dimensional and real, and Lie group homomorphisms are smooth.
 
\subsection{Graded algebras and horizontal homomorphisms}
 
A vector space $V$ is said to be \emph{graded} if it admits a finite direct sum decomposition
$V=\bigoplus_{j\geq1}V_j$ and we refer to the subspaces $V_i$ as layers. We write
$V_{\leq i}$, $V_{<i}$ and $V_{>i}$ for the corresponding sums of
layers, and $\operatorname{pr}_{V_j}\colon V\to V_j$ for the projection
along the other layers. A subspace $U\subseteq V$ is \emph{homogeneous}
if $U=\bigoplus_j(U\cap V_j)$, equivalently if it contains the
homogeneous components of its elements. A linear map
$A\colon V\to V'$ between graded vector spaces is \emph{graded} if
$A(V_j)\subseteq V'_j$ for every $j$. As a direct consequence, the kernel and image of a graded linear map are also
homogeneous.
 
A \emph{graded Lie algebra}  is a graded vector space $V$ with a Lie bracket
satisfying $[V_i,V_j]\subseteq V_{i+j}$. It is a \emph{Carnot algebra}
if, moreover, $V_{j+1}=[V_1,V_j]$ for every $j\geq1$, that is, if it is
generated by $V_1$. Such a grading is called a \textit{stratification} of the Lie algebra. The \emph{step of a grading} is the largest $j$ with
$V_j\neq\{0\}$ and it coincides with the nilpotency step of the Lie algebra in the case of a stratification. In a graded Lie algebra, the ideal generated by a set
of homogeneous elements is homogeneous. If $I$ is a homogeneous ideal,
then $V/I$ is graded by the images of the layers $V_j$, and it is a
Carnot algebra whenever $V$ is.

We now return to the increasing filtration $(\mathcal{F}_iE)_{i\ge 0}$ defined by the subspace $\Delta$ of the vector space $E$ underlying the Lie algebra $\mathfrak g=(E,[\cdot,\cdot]_\mathfrak g)$.
The inclusion
$[\mathcal F_iE,\mathcal F_jE]_{\mathfrak g}\subseteq\mathcal F_{i+j}E$
is immediate from the definition, since the bracket of two bracket words
of lengths $i$ and $j$ is a bracket word of length $i+j$. Hence
$[\cdot,\cdot]_{\mathfrak g}$ induces a graded bracket on $\mathfrak v:=\operatorname{gr}_\Delta\mathfrak g$.
By the Jacobi identity, every bracket word of length $j$ is a linear
combination of brackets $[x_1,[x_2,[\dots,x_j]\dots]]$ with
$x_i\in\Delta$, so $\mathfrak v_j=[\mathfrak v_1,\mathfrak v_{j-1}]$ and
$\mathfrak v$ is a Carnot algebra of step $s$.
 
We say a Lie algebra homomorphism $L\colon\mathfrak g\to\mathfrak h$ into a
graded Lie algebra $(\mathfrak h,[\cdot,\cdot]_\mathfrak h)$ is \emph{horizontal} if
$L(\Delta)\subseteq\mathfrak h_1$. Notice that in this definition, in
Lemmata~\ref{lem:horizontal-gr}, \ref{lem:gr-adapted} and ~\ref{lem:kernel-criterion},
only the grading of $\mathfrak h$ is used, so $\mathfrak h$ may be any
graded Lie algebra, not necessarily Carnot. These results will also be applied to horizontal maps from $\mathfrak g$ into $\mathfrak G$ instead of
$\mathfrak h$. Elsewhere, the target Lie algebra is required to be the Lie algebra of a Carnot
group, and the grading considered is necessarily the stratification.

\begin{lemma}\label{lem:horizontal-gr}
Let $L\colon\mathfrak g\to\mathfrak h$ be a horizontal Lie algebra
homomorphism into a graded Lie algebra.
\begin{enumerate}
 \item[(i)] $L(\mathcal F_iE)\subseteq\mathfrak h_{\leq i}$ for every $i$.
 \item[(ii)] The formula
       $(\operatorname{gr}L)(x+\mathcal F_{i-1}E)
        :=\operatorname{pr}_{\mathfrak h_i}L(x)$ for $x\in\mathcal F_iE$
       defines a graded Lie algebra homomorphism
       $\operatorname{gr}L\colon\mathfrak v\to\mathfrak h$.
 \item[(iii)] $L|_\Delta=(\operatorname{gr}L)|_{\mathfrak v_1}$, and
       two horizontal homomorphisms with the same associated graded map
       coincide.
\end{enumerate}
\end{lemma}
 
\begin{proof}
(i) The map $L$ sends a bracket word of length $j\leq i$ in $\Delta$ to
a bracket of $j$ elements of $\mathfrak h_1$, which lies in
$\mathfrak h_j\subseteq\mathfrak h_{\leq i}$.
 
(ii) By (i), $L(\mathcal F_{i-1}E)\subseteq\mathfrak h_{\leq i-1}$ has
no component in $\mathfrak h_i$, so $\operatorname{gr}L$ is well
defined, and it is graded by construction. Let $x\in\mathcal F_iE$ and
$y\in\mathcal F_jE$. By (i), $L(x)\in\mathfrak h_{\leq i}$ and
$L(y)\in\mathfrak h_{\leq j}$, so the component of weight $i+j$ of
$L([x,y]_{\mathfrak g})=[L(x),L(y)]_{\mathfrak h}$ is
$[\operatorname{pr}_{\mathfrak h_i}L(x),\operatorname{pr}_{\mathfrak h_j}L(y)]_{\mathfrak h}$.
Since the bracket of the classes of $x$ and $y$ in $\mathfrak v$ is the
class of $[x,y]_{\mathfrak g}$ in $\mathfrak v_{i+j}$, this shows that
$\operatorname{gr}L$ preserves brackets.
 
(iii) For $x\in\Delta=\mathcal F_1E$ we have
$(\operatorname{gr}L)(x)=\operatorname{pr}_{\mathfrak h_1}L(x)=L(x)$.
Two Lie algebra homomorphisms that agree on $\Delta$ agree on the Lie
algebra it generates, which is $\mathfrak g$.

\end{proof}
 
\subsection{Adapted gradings and osculating algebras}

\begin{definition}\label{def:adapted-grading}
Given an increasing filtration $(\mathcal{F}_iE)_{i\ge 0}$ of a vector space $E$, an \emph{adapted grading} is a grading $E=\bigoplus_{j=1}^sW_j$ such that
$\mathcal F_jE=\bigoplus_{i=1}^{j}W_i$ for every $j$. In particular
$W_1=\mathcal{F}_1E=\Delta$.
\end{definition}

We do not require $[W_i,W_j]_{\mathfrak g}\subseteq W_{i+j}$. Given a filtration $(\mathcal{F}_iE)_{i\ge 0}$, it is always possible to construct an adapted
grading: take $W_1=\mathcal{F}_1E=\Delta$ and any complement
$W_j$ of $\mathcal F_{j-1}E$ in $\mathcal F_jE$ for $j\geq2$. 

Since
$\mathcal F_jE=\mathcal F_{j-1}E\oplus W_j$, every
$\xi\in\mathfrak v_j$ has a unique representative in $W_j$, and we write
$\sigma_W\colon\mathfrak v\to E$ for the linear isomorphism that sends
$\xi$ to this representative. Transporting the bracket of $\mathfrak v$ by $\sigma_W$
gives the \emph{osculating algebra}
$\mathfrak g_0^W=(E,[\cdot,\cdot]_{0,W})$, where
$[u,v]_{0,W}:=\sigma_W\bigl([\sigma_W^{-1}(u),\sigma_W^{-1}(v)]_{\mathfrak v}\bigr)$ for any $u,v\in E$.
We grade $\mathfrak g_0^W$ by the subspaces $W_j$, and define the corresponding homogeneous
dilations
\[
 \delta^W_t\Bigl(\sum_{j=1}^sw_j\Bigr):=\sum_{j=1}^st^jw_j\ ,
 \quad w_j\in W_j,\ t>0\,.
\]

As a way to measure the difference between the original bracket $[\cdot,\cdot]_\mathfrak g$ and the osculating bracket $[\cdot,\cdot]_{0,W}$, we introduce the \emph{bracket defect} of $W$ as the bilinear map
\[
 B_W\colon E\times E\to E\ ,\quad
 B_W(u,v):=[u,v]_{\mathfrak g}-[u,v]_{0,W}\,.
\]
\begin{lemma}\label{lem:gr-adapted}
Let $L\colon\mathfrak g\to\mathfrak h$ be a Lie algebra homomorphism
into a graded Lie algebra, and let $W$ be an adapted grading with
$L(W_j)\subseteq\mathfrak h_j$ for every $j$. Then $L$ is horizontal
and $\operatorname{gr}L=L\circ\sigma_W$.
\end{lemma}

\begin{proof}
$L$ is horizontal because, by construction, $W_1=\Delta$. For $\xi\in\mathfrak v_j$, the
vector $\sigma_W(\xi)\in W_j\subseteq\mathcal F_jE$ represents $\xi$,
and $L(\sigma_W(\xi))\in\mathfrak h_j$. Hence
$(\operatorname{gr}L)(\xi)=\operatorname{pr}_{\mathfrak h_j}L(\sigma_W(\xi))=L(\sigma_W(\xi))$.

\end{proof}
\begin{lemma}\label{lem:osculating}
Let $W$ be an adapted grading.
\begin{enumerate}
 \item[(i)] The map $\sigma_W\colon\mathfrak v\to\mathfrak g_0^W$ is a graded
       Lie algebra isomorphism. Hence $\mathfrak g_0^W$ is a Carnot
       algebra of step $s$ with first layer $\Delta$, and the dilations
       $\delta^W_t$ are automorphisms of $\mathfrak g_0^W$.
 \item[(ii)] For $u\in W_i$ and $v\in W_j$,
       \[
        [u,v]_{0,W}=\operatorname{pr}_{W_{i+j}}[u,v]_{\mathfrak g}
        \ \text{and}\ 
        B_W(u,v)\in\mathcal F_{i+j-1}E\,,
       \]
       with the conventions $W_k=\{0\}$ and $\mathcal F_kE=E$ for
       $k>s$.
 \item[(iii)] For all $u,v\in E$,
       $\delta^W_{1/t}[\delta^W_tu,\delta^W_tv]_{\mathfrak g}$ is
       polynomial in $t$ and converges to $[u,v]_{0,W}$ as
       $t\to 0$.
\end{enumerate}
\end{lemma}
 
\begin{proof}
(i) This holds by the definition of $[\cdot,\cdot]_{0,W}$, since
$\sigma_W(\mathfrak v_j)=W_j$.
 
(ii) The bracket in $\mathfrak v$ of the classes of $u$ and $v$ is the
class of $[u,v]_{\mathfrak g}\in\mathcal F_{i+j}E=W_{\leq i+j}$ in
$\mathfrak v_{i+j}$, and its representative in $W_{i+j}$ is
$\operatorname{pr}_{W_{i+j}}[u,v]_{\mathfrak g}$. If $i+j>s$, both sides
vanish. In either case, $B_W(u,v)$ is the sum of the components of
$[u,v]_{\mathfrak g}$ of weight less than $i+j$.
 
(iii) By bilinearity, it suffices to take $u\in W_i$ and $v\in W_j$, for which
\[
 \delta^W_{1/t}[\delta^W_tu,\delta^W_tv]_{\mathfrak g}
 =\sum_{k\leq i+j}t^{\,i+j-k}\operatorname{pr}_{W_k}[u,v]_{\mathfrak g}\,,
\]
which converges to $\operatorname{pr}_{W_{i+j}}[u,v]_{\mathfrak g}$ as $t\to 0$. The claim then follows since $\operatorname{pr}_{W_{i+j}}[u,v]_\mathfrak g=[u,v]_{0,W}$ by (ii).

\end{proof}
 
Part (iii) shows that $\mathfrak g_0^W$ is the limit of the brackets
rescaled by $\delta^W_t$. Its simply connected Lie group, the
\emph{osculating group} \cite{FS74,vEY17}, is the metric tangent cone of
$G$ by Mitchell's theorem (Theorem~\ref{thm:mitchell} below).
 
\begin{lemma}\label{lem:grading-change-algebra}
Let $W$ and $W'$ be adapted gradings and define
$C:=\sigma_{W'}\circ\sigma_W^{-1}\colon E\to E$. Then $C(W_j)=W'_j$ for
every $j$, $C$ is a graded Lie algebra isomorphism
$\mathfrak g_0^W\to\mathfrak g_0^{W'}$,
$C\circ\delta^W_t=\delta^{W'}_t\circ C$ for every $t>0$, and
$(C-\mathrm{id})(W_j)\subseteq\mathcal F_{j-1}E$.
\end{lemma}
 
\begin{proof}
The first three assertions follow from Lemma~\ref{lem:osculating}(i)
for $W$ and $W'$. For $u\in W_j$, the vectors $u$ and $C(u)$ represent
the same class $\sigma_W^{-1}(u)\in\mathfrak v_j$, so their difference
lies in $\mathcal F_{j-1}E$.

\end{proof}
 
\subsection{Defect ideals}
\begin{definition}\label{def:defect-ideal}
Let $W$ be an adapted grading. The \emph{defect ideal} of $W$ is the
ideal of $\mathfrak g_0^W$ generated by the homogeneous components of
the bracket defects,
\[
 J_W:=\operatorname{Ideal}_{0,W}
 \bigl(\{\operatorname{pr}_{W_i}B_W(u,v):
 u,v\in E,\ 1\leq i\leq s\}\bigr)\,,
\]
and we define $j_W:=\sigma_W^{-1}(J_W)\subseteq\mathfrak v$.
\end{definition}

Being generated by homogeneous elements, $J_W$ is a homogeneous ideal of
$\mathfrak g_0^W$. It is therefore the smallest subspace of $E$ that
contains $B_W(u,v)$ for all $u,v\in E$, is closed under
$[E,\cdot\,]_{0,W}$, and is homogeneous with respect to $W$. Since
$\sigma_W$ is a graded isomorphism, $j_W$ is a homogeneous ideal of
$\mathfrak v$.
 
\begin{lemma}[Kernel criterion]\label{lem:kernel-criterion}
Let $D\colon\mathfrak v\to\mathfrak h$ be a graded Lie algebra
homomorphism into a graded Lie algebra, let $W$ be an adapted grading,
and put $A_W:=D\circ\sigma_W^{-1}$. Then
\[
 A_W\text{ preserves }[\cdot,\cdot]_{\mathfrak g}
 \iff J_W\subseteq\ker A_W
 \iff j_W\subseteq\ker D\,.
\]
In that case $A_W$ is a horizontal Lie algebra homomorphism with
$\operatorname{gr}A_W=D$.
\end{lemma}
 
\begin{proof}
By (i) in Lemma~\ref{lem:osculating}, $A_W\colon \mathfrak g_0^W\to\mathfrak h$ is a graded homomorphism. Hence, for $u,v\in E$,
\[
 A_W([u,v]_{\mathfrak g})-[A_W(u),A_W(v)]_{\mathfrak h}
 =A_W([u,v]_{\mathfrak g})-A_W([u,v]_{0,W})
 =A_W(B_W(u,v))\,,
\]
so $A_W$ preserves $[\cdot,\cdot]_{\mathfrak g}$ if and only if it
annihilates every bracket defect. The kernel of $A_W$ is a homogeneous
ideal of $\mathfrak g_0^W$, so this happens if and only if
$J_W\subseteq\ker A_W$. Applying $\sigma_W^{-1}$ gives the second
equivalence. Finally, $A_W(W_j)=D(\mathfrak v_j)\subseteq\mathfrak h_j$,
so $A_W$ is horizontal, and
$\operatorname{gr}A_W=A_W\circ\sigma_W=D$ by
Lemma~\ref{lem:gr-adapted}.

\end{proof}
 
\subsection{Groups, exponential coordinates and distances}

For a nilpotent Lie algebra $\mathfrak a$ with Lie bracket $[\cdot,\cdot]_\frak a$, we write $\star_{\mathfrak a}$
for the Baker-Campbell-Hausdorff (BCH) product
\[
 u\star_{\mathfrak a}v=u+v+\tfrac12[u,v]_{\mathfrak a}
 +\tfrac1{12}\bigl([u,[u,v]_{\mathfrak a}]_{\mathfrak a}
 +[v,[v,u]_{\mathfrak a}]_{\mathfrak a}\bigr)+\cdots\,,
\]
which is a finite sum. Then $(\mathfrak a,\star_{\mathfrak a})$ is the
simply connected Lie group with Lie algebra $\mathfrak a$, in
exponential coordinates. Its identity element is $0$ and the inverse of
$u$ is $-u$.

Since the group $G$ considered is simply connected and nilpotent, the
exponential map $\mathfrak g\to G$ is a diffeomorphism, and so we can
identify $G$ with $(E,\star_{\mathfrak g})$, and the Haar measure is chosen to be the
Lebesgue measure on $E$. In the same way, we have the Carnot group $H=(\mathfrak h,\star_{\mathfrak h})$,
and the \emph{osculating group} is $G_0^W:=(E,\star_{0,W})$, where
$\star_{0,W}$ is the BCH product of $[\cdot,\cdot]_{0,W}$. The dilations
$\delta^W_t$ are automorphisms of $G_0^W$, and the dilations
$\bar\delta_t(\sum_jh_j)=\sum_jt^jh_j$, with $h_j\in\mathfrak h_j$, are
automorphisms of $H$.

\begin{lemma}\label{lem:linear-homs}
Let $\mathfrak a$ and $\mathfrak b$ be nilpotent Lie algebras. A map
$\Phi\colon\mathfrak a\to\mathfrak b$ is a Lie group homomorphism
$(\mathfrak a,\star_{\mathfrak a})\to(\mathfrak b,\star_{\mathfrak b})$
if and only if it is a Lie algebra homomorphism. In particular, a linear
map $A\colon E\to\mathfrak h$ satisfies
$A(u\star_{\mathfrak g}v)=A(u)\star_{\mathfrak h}A(v)$ for all
$u,v\in E$ if and only if it preserves $[\cdot,\cdot]_{\mathfrak g}$.
\end{lemma}

\begin{proof}
A Lie group homomorphism satisfies
$\Phi\circ\exp=\exp\circ\,d\Phi_0$, and in exponential coordinates both
exponential maps are the identity, so $\Phi=d\Phi_0$, a Lie algebra
homomorphism. Conversely, a Lie algebra homomorphism preserves every
iterated bracket in the BCH formula.

\end{proof}

Fix an inner product on $\Delta$ and extend it to the left-invariant
distribution $x\mapsto d(L_x)_0(\Delta)$, where
$L_x(y)=x\star_{\mathfrak g}y$. The
Carnot-Carath\'eodory distance $d_G$ is the infimum of the lengths of
horizontal curves, and by the Chow-Rashevskii theorem it is finite,
left-invariant, and induces the manifold topology. 
In the same way,
$d_{0,W}$ denotes the distance on $G_0^W$ defined by $W_1=\Delta$ with
the same inner product, and $d_H$ the distance on $H$ defined by an
inner product on $\mathfrak h_1$. Since $\delta^W_t$ and $\bar\delta_t$
are automorphisms that multiply the first layer by $t$,
\[
 d_{0,W}(\delta^W_tu,\delta^W_tv)=t\,d_{0,W}(u,v)\ \text{ and }\ d_H(\bar\delta_ta,\bar\delta_tb)=t\,d_H(a,b)\,.
\]

Other choices of inner products give bi-Lipschitz equivalent distances,
so the class of Lipschitz maps $G\to H$ and the notion of Pansu
differentiability below do not depend on them.

For $t>0$ put
$d^W_t(u,v):=\tfrac1t\,d_G\bigl(\delta^W_tu,\delta^W_tv\bigr)$. The
dilations $\delta^W_t$ need not be automorphisms of $G$, so $d^W_t$ will not coincide with $d_G$ in general, however the following holds.
\begin{theorem}[Mitchell's Theorem, see {\cite[Theorem~6.1]{LDNGT}}]
\label{thm:mitchell}
For every adapted grading $W$ and every compact set $K\subseteq E$,
\begin{equation}\label{eq:dilated-convergence}
 \lim_{t\to 0}\ \sup_{u,v\in K}
 \bigl|d^W_t(u,v)-d_{0,W}(u,v)\bigr|=0\,.
\end{equation}
\end{theorem}

\begin{proof}
For $u,v$ in a neighbourhood $U$ of $0$, this is
\cite[Theorem~6.1]{LDNGT}, which also gives a rate of convergence. For an arbitrary compact set $K$, choose $\lambda>0$ with
$\delta^W_\lambda(K)\subseteq U$. Since
$\delta^W_{t/\lambda}\circ\delta^W_\lambda=\delta^W_t$, the definition of
$d^W_t$ gives
\[
 d^W_{t/\lambda}\bigl(\delta^W_\lambda u,\delta^W_\lambda v\bigr)
 =\frac{\lambda}{t}\,d_G\bigl(\delta^W_tu,\delta^W_tv\bigr)
 =\lambda\,d^W_t(u,v)\,.
\]
By the homogeneity of $d_{0,W}$, we then get
\[
 d^W_t(u,v)-d_{0,W}(u,v)
 =\lambda^{-1}\Bigl(d^W_{t/\lambda}\bigl(\delta^W_\lambda u,\delta^W_\lambda v\bigr)
 -d_{0,W}\bigl(\delta^W_\lambda u,\delta^W_\lambda v\bigr)\Bigr)\ \text{ for all }u,v\in K\,.
\]
Since $\delta^W_\lambda u,\delta^W_\lambda v\in U$, the right-hand side tends to
$0$ as $t\to 0$ uniformly on $K\times K$ by the first case.

\end{proof}
\begin{corollary}\label{cor:ball-box}
For every adapted grading $W$ there is $r>0$ such that
$\tfrac12\,d_{0,W}(0,u)\leq d_G(0,u)\leq2\,d_{0,W}(0,u)$ whenever
$d_{0,W}(0,u)\leq r$.
\end{corollary}

\begin{proof}
The unit sphere $\Sigma$ of $d_{0,W}$ about $0$ is compact, since
$d_{0,W}$ is homogeneous and induces the manifold topology.
Theorem~\ref{thm:mitchell}, applied to $K=\Sigma\cup\{0\}$, gives $r>0$
with $\tfrac12\leq d^W_t(0,v)\leq\tfrac32$ for all $v\in\Sigma$ and
$0<t\leq r$. If $0<t:=d_{0,W}(0,u)\leq r$, then
$v:=\delta^W_{1/t}u\in\Sigma$ by homogeneity, and so
\[
 d_G(0,u)=d_G\bigl(\delta^W_t0,\delta^W_tv\bigr)=t\,d^W_t(0,v)
 \in\bigl[\tfrac t2,\tfrac{3t}2\bigr]\,,
\]
which yields the claim.

\end{proof}
\subsection{Pansu differentiability}\label{subsec: pansu diff}

\begin{definition}\label{def:pansu-W}
Let $\Omega\subseteq G$ be open, $f\colon\Omega\to H$, $x\in\Omega$, and
let $W$ be an adapted grading. For $t>0$, consider
\[
 R^W_{x,t}(u):=\bar\delta_{1/t}\bigl(f(x)^{-1}\star_{\mathfrak h}
 f(x\star_{\mathfrak g}\delta^W_tu)\bigr)\,,
\]
which is defined as long as $x\star_{\mathfrak g}\delta^W_tu\in\Omega$.
Since $\Omega$ is open and $\delta^W_t\to0$ uniformly on compact sets
as $t\to 0$,  this holds for all
$u$ in any compact $K\subseteq E$ once $t$ is small enough.

We say the map $f$ is \emph{Pansu differentiable at $x$ with respect to $W$} if
there is a graded Lie algebra homomorphism
$A\colon\mathfrak g_0^W\to\mathfrak h$ such that
\[
 \lim_{t\to 0}\ \sup_{u\in K}
 d_H\bigl(R^W_{x,t}(u),A(u)\bigr)=0
 \quad\text{for every compact set  }K\subseteq E\,.
\]
If this is the case, we write $D^W_Pf(x):=A$. By Lemma~\ref{lem:linear-homs}, it is
also a homogeneous group homomorphism $D^W_Pf(x)\colon G_0^W\to H$.
\end{definition}

If $\mathfrak g$ is itself a Carnot algebra with first layer $\Delta$
and $W$ is its stratification, then $\mathfrak g_0^W=\mathfrak g$, the
dilations $\delta^W_t$ are automorphisms of $G$, and
Definition~\ref{def:pansu-W} is Pansu's definition \cite{Pan89}. In
that case Theorem~\ref{thm:existence} below is Pansu's theorem (a
detailed proof can be found in \cite[Theorem~2.1.6]{Monti2001}). Pansu's theorem
has also been extended in the direction of the target, to graded
groups \cite[Theorem~2.11]{Mag10} and to Banach homogeneous groups
\cite[Theorem~1.1]{MR14}, while Theorem~\ref{thm:existence} below instead relaxes
the assumption on the domain.

\begin{theorem}[Pansu \cite{Pan89} and Karmanova-Vodopyanov \cite{KV08}]\label{thm:existence}
Let $\Omega\subseteq G$ be open, $f\colon\Omega\to H$ Lipschitz, and $W$
an adapted grading. Then $f$ is Pansu differentiable with respect to
$W$ at almost every $x\in\Omega$.
\end{theorem}

\begin{proof}We apply \cite[Theorem~3.3.5]{KV08} to $G$ with a left-invariant
frame adapted to $W$, and explain how the objects used there read in
our setting.  For each $j$, choose a basis of $W_j$ of weight $j$ (referred to as the degree in \cite{KV08}).
Together these form a basis $X_1,\dots,X_n$ of $E$, which is
a frame of left-invariant vector fields on $G$. The flow of a left-invariant vector field is right
multiplication, so the normal coordinates of \cite{KV08} at $x$ are
$y\mapsto x\star_{\mathfrak g}y$, and their dilations are
$x\star_{\mathfrak g}y\mapsto x\star_{\mathfrak g}\delta^W_ty$. The
structure constants of a left-invariant frame are constant, so the
nilpotent tangent cone of \cite[Theorem~2.1.6]{KV08} is simply
$(E,\star_{0,W})$. By \cite[Definition~2.1.18 and Lemma~2.1.22]{KV08},
the local Carnot group at $x$ is $(E,\star_{0,W})$ transported by
$y\mapsto x\star_{\mathfrak g}y$. For the stratified frame of $H$, the
local Carnot group at $f(x)$ is $H$, with distance $d_H$.

A horizontal homomorphism between these local Carnot groups is
therefore $x\star_{\mathfrak g}y\mapsto f(x)\star_{\mathfrak h}A(y)$,
with $A\colon\mathfrak g_0^W\to\mathfrak h$ graded. With these
identifications, hc-differentiability at $x$ \cite[(1.0.1)]{KV08} reads
\[
 d_H\bigl(f(x\star_{\mathfrak g}y),f(x)\star_{\mathfrak h}A(y)\bigr)
 =o\bigl(d_{0,W}(0,y)\bigr)
 \ \text{ as }\ y\to0\,,
\]
where we used Corollary~\ref{cor:ball-box} to replace $d_G(0,y)$ by
$d_{0,W}(0,y)$. Writing $y=\delta^W_tv$ with $d_{0,W}(0,v)=1$, this is
uniform convergence of $R^W_{x,t}$ to $A$ on the unit sphere of
$d_{0,W}$, and, by homogeneity, on every compact set. 

Finally, by the
ball-box theorem a map that is Lipschitz for $d_G$ and $d_H$ is
Lipschitz, on compact sets, for the quasimetrics $d_\infty$ used in
\cite{KV08}. The claim is now obtained from \cite[Theorem~3.3.5]{KV08}.

\end{proof}

We write $\Omega_f^W$ for the set of points at which $f$ is Pansu
differentiable with respect to $W$ which, by Theorem~\ref{thm:existence}
has full measure when $f$ is Lipschitz. Crucially, in the case of a  Lipschitz map $f$, this set and
the differential do not depend on $W$, up to the identification
$\sigma_W$.

\begin{lemma}[Change of grading]\label{lem:change-grading}
Let $f\colon\Omega\subseteq G\to H$ be Lipschitz, let $x\in\Omega$, and let $W$ and
$W'$ be adapted gradings. If $f$ is Pansu differentiable at $x$ with
respect to $W$, then it is Pansu differentiable at $x$ with respect to
$W'$, and
\[
 D^{W'}_Pf(x)\circ\sigma_{W'}=D^W_Pf(x)\circ\sigma_W\,.
\]
\end{lemma}

\begin{proof}
Put $C:=\sigma_{W'}\circ\sigma_W^{-1}$ and $P:=C-\mathrm{id}_E$. By
Lemma~\ref{lem:grading-change-algebra}, $C$ is a graded Lie algebra
isomorphism $\mathfrak g_0^W\to\mathfrak g_0^{W'}$,
$C\circ\delta^W_t=\delta^{W'}_t\circ C$, and
$P(W_j)\subseteq\mathcal F_{j-1}E=W_{\leq j-1}$. Hence
$C\delta^W_tv=\delta^W_t\bigl(v+e_t(v)\bigr)$ for $v\in E$, where
\[
 e_t(v):=\delta^W_{1/t}P\delta^W_tv
 =\sum_j\sum_{i<j}t^{\,j-i}\operatorname{pr}_{W_i}P(v_j)\ ,
 \quad v=\sum_jv_j,\ v_j\in W_j\,,
\]
so $e_t\to0$ uniformly on compact sets as $t\downarrow0$.

Fix a compact set $K\subseteq E$ and a compact neighbourhood $K'$ of
$K$. By left invariance of $d_G$, for $v\in K$ and $t$ small, so that $v,v+e_t(v)\in K'$
\[
 \tfrac1t\,d_G\bigl(x\star_{\mathfrak g}C\delta^W_tv,\,
 x\star_{\mathfrak g}\delta^W_tv\bigr)
 =d^W_t\bigl(v+e_t(v),v\bigr)
 \leq d_{0,W}\bigl(v+e_t(v),v\bigr)
 +\sup_{K'\times K'}\bigl|d^W_t-d_{0,W}\bigr|\,,
\]
which tends to $0$ uniformly in $v\in K$, by
Theorem~\ref{thm:mitchell} and the uniform continuity of $d_{0,W}$ on
$K'$. Since $\delta^{W'}_tCv=C\delta^W_tv$, left invariance and
homogeneity of $d_H$ and the Lipschitz bound for $f$ give
\[
 d_H\bigl(R^{W'}_{x,t}(Cv),R^W_{x,t}(v)\bigr)
 =\tfrac1t\,d_H\bigl(f(x\star_{\mathfrak g}C\delta^W_tv),
 f(x\star_{\mathfrak g}\delta^W_tv)\bigr)
 \longrightarrow0
\]
uniformly in $v\in K$. Hence $R^{W'}_{x,t}\circ C$ converges to
$D^W_Pf(x)$ uniformly on compact sets, that is, $R^{W'}_{x,t}$ converges
to $D^W_Pf(x)\circ C^{-1}$ uniformly on compact sets. The latter is a
graded Lie algebra homomorphism $\mathfrak g_0^{W'}\to\mathfrak h$, so
$D^{W'}_Pf(x)=D^W_Pf(x)\circ C^{-1}$, which is the claim.
\end{proof}

\begin{definition}\label{def:canonical-differential}
Let $f\colon\Omega\to H$ be Lipschitz. By
Lemma~\ref{lem:change-grading}, the set
$\Omega_f:=\Omega^W_f$ does not depend on $W$, and it has full measure.
For $x\in\Omega_f$, the \emph{canonical Pansu differential} of $f$ at
$x$ is the graded Lie algebra homomorphism
\[
 D_Pf(x):=D^W_Pf(x)\circ\sigma_W\colon\mathfrak v\to\mathfrak h\,,
\]
which does not depend on $W$. Equivalently,
$D^W_Pf(x)=D_Pf(x)\circ\sigma_W^{-1}$ for every adapted grading $W$.
\end{definition}

\begin{example}\label{ex:heisenberg-differentials}
Let $\mathfrak g$, $W$ and $\widetilde W$ be as in the tilted Heisenberg
example of the introduction, and write $u=xX+yY+zZ$. Then
$C:=\sigma_{\widetilde W}\circ\sigma_W^{-1}$ fixes $X$ and $Y$ and maps
$Z$ to $T=Z+X$, and
\[
 \delta^{\widetilde W}_tu=\bigl(t(x-z)+t^2z\bigr)X+tyY+t^2zZ\,.
\]
Let the target group be again the first Heisenberg group $H=G$ with the stratification $W$ and let $f\colon G\to G$ be the identity map. Since
$f(p)^{-1}\star_{\mathfrak g}f(p\star_{\mathfrak g}v)=v$ at every
$p\in G$, we get $R^W_{p,t}(u)=u$ and
\[
 R^{\widetilde W}_{p,t}(u)=\bar\delta_{1/t}\bigl(\delta^{\widetilde W}_tu\bigr)
 =(x-z+tz)X+yY+zZ
 \longrightarrow(x-z)X+yY+zZ\,.
\]
Thus $D^W_P\mathrm{id}(p)=\mathrm{id}_E$. On the other hand,
$D^{\widetilde W}_P\mathrm{id}(p)=C^{-1}$, which maps $X,Y,T$ to
$X,Y,Z$, in accordance with Lemma~\ref{lem:change-grading}. Both represent the same
canonical differential $D_P\mathrm{id}(p)=\sigma_W$. The map $C^{-1}$ is
a graded Lie algebra homomorphism $\mathfrak g_0^{\widetilde W}\to\mathfrak g$,
but it does not preserve $[\cdot,\cdot]_{\mathfrak g}$: it maps
$[X,Y]_{\mathfrak g}=Z$ to $Z-X$. So the $\widetilde W$-differential of
the identity is not a group homomorphism $G\to H$ at any point. This is
consistent with Theorem~\ref{thm:A}, since $W$ is minimal and
$\widetilde W$ is not.
\end{example}

\section{Two quotients}\label{sec:quotients}

This section constructs the universal horizontal quotient
$(\mathfrak G,\rho)$ of the introduction, which depends only on
$(\mathfrak g,\Delta)$, and, for each adapted grading $W$, the common
quotient $\mathfrak q_W:=E/J_W$ of $\mathfrak g$ and $\mathfrak g_0^W$.
Their universal properties concern different classes of
homomorphisms: horizontal ones for $\mathfrak G$, and those that
preserve the layers of $W$ for $\mathfrak q_W$.
Section~\ref{subsec:comparison} shows that $\mathfrak q_W$ is a quotient
of $\mathfrak G$ and identifies the kernel.

\subsection{The universal horizontal quotient}

Let $\mathcal L(\Delta)=\bigoplus_{j\geq1}\mathcal L_j(\Delta)$ be the
free Lie algebra on the vector space $\Delta$, graded by bracket length,
so that $\mathcal L_1(\Delta)=\Delta$ and $\mathcal L_j(\Delta)$ is
spanned by the brackets of $j$ elements of $\Delta$. According to our definition, these are all elements of weight $j$. Every linear map
from $\Delta$ to a Lie algebra extends uniquely to a Lie algebra
homomorphism on $\mathcal L(\Delta)$ \cite{Reu93}. Let
\[
 \operatorname{ev}\colon\mathcal L(\Delta)\to\mathfrak g
\]
be the homomorphism extending the inclusion $\Delta\subseteq E$: it
evaluates formal brackets of elements of $\Delta$ in $\mathfrak g$. Since
$\mathcal F_iE$ is spanned by the bracket words of length at most $i$
in $\Delta$, we have
\begin{equation}\label{eq:epsilon-filtration}
 \operatorname{ev}\left(\bigoplus_{j\leq i}\mathcal L_j(\Delta)\right)
 =\mathcal F_iE
 \quad\text{for every }i\,.
\end{equation}
In particular, the homomorphism $\operatorname{ev}$ is surjective. 

Given $R:=\ker\operatorname{ev}$, let $R^{\mathrm h}$ be the ideal of $\mathcal L(\Delta)$ generated
by the homogeneous components of the elements of $R$. Define
\[
 \mathfrak G:=\mathcal L(\Delta)/R^{\mathrm h}\ \text{ and }
 \ 
 \rho\colon\mathfrak g\to\mathfrak G\ ,
 \quad
 \rho(\operatorname{ev}(\xi)):=\xi+R^{\mathrm h}\,.
\]
The map $\rho\colon\mathfrak g\to\mathfrak G$ is well-defined because $R\subseteq R^{\mathrm h}$. Moreover, we can equip  $\mathfrak G$ with a grading by the images $\mathfrak G_i$ of the layers
$\mathcal L_i(\Delta)$.

\begin{proposition}[Universal horizontal quotient]
\label{prop:horizontal-quotient}
The algebra $\mathfrak G$ is a Carnot algebra with $\mathfrak G_j=\{0\}$
for any $j>s$ where $s$ is the distribution step, and $\rho\colon\mathfrak g\to\mathfrak G$ is a surjective Lie algebra homomorphism with
$\rho(\Delta)=\mathfrak G_1$. Moreover:
\begin{enumerate}
 \item[(i)] The map $\rho$ is strict:
       \begin{equation}\label{eq:rho-strict}
        \rho(\mathcal F_iE)=\mathfrak G_{\leq i}
        \quad\text{for every }i\,.
       \end{equation}
       Hence $\operatorname{gr}\rho\colon\mathfrak v\to\mathfrak G$ is a
       graded Lie algebra homomorphism mapping $\mathfrak v_j$ onto
       $\mathfrak G_j$ for every $j$.
 \item[(ii)] For every Carnot algebra $\mathfrak h$ and every horizontal Lie
       algebra homomorphism $L\colon\mathfrak g\to\mathfrak h$, there is
       a unique graded Lie algebra homomorphism
       $\overline L\colon\mathfrak G\to\mathfrak h$ with
       $L=\overline L\circ\rho$. It satisfies
       $\operatorname{gr}L=\overline L\circ\operatorname{gr}\rho$.
\end{enumerate}
\end{proposition}

\begin{proof}
If $j>s$ and $\xi\in\mathcal L_j(\Delta)$, then
$\operatorname{ev}(\xi)\in E=\mathcal F_sE$, so by
\eqref{eq:epsilon-filtration} there is
$\eta\in\bigoplus_{k\leq s}\mathcal L_k(\Delta)$ with
$\operatorname{ev}(\eta)=\operatorname{ev}(\xi)$. Thus $\xi-\eta\in R$,
and its component of weight $j$ is $\xi$, since $\eta$ has no component
of weight $j>s$. Hence $\xi\in R^{\mathrm h}$, and
$\mathfrak G_j=\{0\}$ for $j>s$.

The
grading of $\mathfrak G$ is a Lie algebra grading, and $\mathfrak G$ is
generated by $\mathfrak G_1$ because $\mathcal L(\Delta)$ is generated by
$\Delta$, hence $\mathfrak G$ is a Carnot algebra. Furthermore, since $\operatorname{ev}$ is
a surjective homomorphism and $\rho\circ\operatorname{ev}$ is the quotient map
of $\mathcal L(\Delta)$, the map $\rho$ is also a surjective homomorphism.

(i) The inclusion $\rho(\mathcal{F}_iE)\subseteq \mathfrak G_{\le i}$ follows directly by construction. Conversely, every element of
$\mathfrak G_{\leq i}$ is $\rho(\operatorname{ev}(\eta))$ for
some $\eta\in\bigoplus_{k\leq i}\mathcal L_k(\Delta)$, and $\operatorname{ev}(\eta)\in\mathcal F_iE$. This
proves \eqref{eq:rho-strict}. For $i=1$ it reads
$\rho(\Delta)=\mathfrak G_1$, i.e. $\rho$ is horizontal, and so
Lemma~\ref{lem:horizontal-gr}(ii) implies that $\operatorname{gr}\rho$ is a graded Lie algebra
homomorphism. Finally, given $\zeta\in\mathfrak G_j$, we can choose $x\in\mathcal F_jE$ with
$\rho(x)=\zeta$, and so
$\operatorname{gr}\rho(x+\mathcal F_{j-1}E)
=\operatorname{pr}_{\mathfrak G_j}\rho(x)=\zeta$, i.e.
$\operatorname{gr}\rho$ maps $\mathfrak v_j$ onto $\mathfrak G_j$.

(ii) The homomorphism $L\circ\operatorname{ev}$ maps $\mathcal L_j(\Delta)$ into
$\mathfrak h_j$, because $\mathcal L_j(\Delta)$ is spanned by brackets of
$j$ elements of $\Delta$ and $L(\Delta)\subseteq\mathfrak h_1$. Moreover, it
vanishes on $R$, and so if $r=\sum_jr_j\in R$ with $r_j\in\mathcal L_j(\Delta)$,
the vectors $L\circ\operatorname{ev}(r_j)\in\mathfrak h_j$ lie in different layers, so each vanishes. Hence the ideal $\ker(L\circ\operatorname{ev})$
contains $R^{\mathrm h}$, and $L\circ\operatorname{ev}$ induces a graded
homomorphism $\overline L\colon\mathfrak G\to\mathfrak h$ with
$\overline L\circ\rho\circ\operatorname{ev}=L\circ\operatorname{ev}$, that is,
$\overline L\circ\rho=L$. It is unique because $\rho$ is surjective.
Finally, for $x\in\mathcal F_jE$, since $\overline L$ is graded,
\[
 \operatorname{gr}L(x+\mathcal F_{j-1}E)
 =\operatorname{pr}_{\mathfrak h_j}\overline L\circ\rho(x)
 =\overline L\circ\operatorname{pr}_{\mathfrak G_j}\rho(x)
 =\overline L\circ \operatorname{gr}\rho(x+\mathcal F_{j-1}E)\,.
\]

\end{proof}

\begin{corollary}\label{cor:horizontal-quotient}$\phantom{=}$
\begin{enumerate}
 \item[(i)] $\ker\rho$ is the intersection of the kernels of all horizontal
       homomorphisms from $\mathfrak g$ into Carnot algebras.
 \item[(ii)] Let $\mathfrak G'$ be a Carnot algebra and
       $\rho'\colon\mathfrak g\to\mathfrak G'$ a horizontal homomorphism
       such that every horizontal homomorphism $L$ from $\mathfrak g$
       into a Carnot algebra is $L=\overline L'\circ\rho'$ for a unique
       graded homomorphism $\overline L'$. Then there is a unique graded
       isomorphism $\phi\colon\mathfrak G\to\mathfrak G'$ with
       $\rho'=\phi\circ\rho$.
 \item[(iii)] If a set $S$ generates $R$ as an ideal, then $R^{\mathrm h}$ is
       the ideal generated by the homogeneous components of the
       elements of $S$.
\end{enumerate}
\end{corollary}

\begin{proof}
(i) Every such $L$ satisfies $\ker\rho\subseteq\ker L$, since
$L=\overline L\circ\rho$, and $\rho$ is itself one of them.

(ii) Proposition~\ref{prop:horizontal-quotient}(ii) applied to $\rho'$
gives a graded homomorphism $\phi$ with $\rho'=\phi\circ\rho$, and by assumption there exists a graded homomorphism $\psi$ such that
$\rho=\psi\circ\rho'$. Then $\psi\circ\phi\circ\rho=\rho$, and so
$\psi\circ\phi=\mathrm{id}$ by the surjectivity of $\rho$. Moreover, by the uniqueness assumption
$\phi\circ\psi\circ\rho'=\rho'$ implies $\phi\circ\psi=\mathrm{id}$. The map $\phi$ is unique because $\rho$ is
surjective.

(iii) Let $I$ be the ideal generated by the homogeneous components of
the elements of $S$. Since $S\subseteq R$, these components lie in
$R^{\mathrm h}$, so $I\subseteq R^{\mathrm h}$. Conversely, $I$ contains
every element of $S$, so $R\subseteq I$ and since $I$ is homogeneous, it
contains the homogeneous components of the elements of $R$, i.e.
$R^{\mathrm h}\subseteq I$.
\end{proof}

Part (iii) is how $\mathfrak G$ is computed in practice: a presentation
of $\mathfrak g$ with generators in $\Delta$ gives a presentation of
$\mathfrak G$ with the same generators, each relation being replaced by
its homogeneous components. For example, take the Heisenberg algebra
of the introduction, but with $\Delta=E$. Then the generators $X,Y,Z$
all have weight $1$, so the relation $[X,Y]-Z$ has the homogeneous
components $[X,Y]$, of weight $2$, and $-Z$, of weight $1$. By (iii),
both vanish in $\mathfrak G$, hence $\mathfrak G$ is abelian of
dimension $2$, spanned by the images of $X$ and $Y$, and $\rho(Z)=0$. In
particular, this shows that $\rho|_\Delta$ need not be injective when
$\Delta\cap[E,E]_{\mathfrak g}\neq\{0\}$.

\begin{definition}\label{def:rstar}
The \emph{visibility weight} of $(\mathfrak g,\Delta)$ is the lowest
weight of the homogeneous ideal $\ker(\operatorname{gr}\rho)$,
\[
 r_*:=\min\bigl\{j\geq1:\ker(\operatorname{gr}\rho)\cap\mathfrak v_j\neq\{0\}\bigr\}\,,
\]
with $r_*=\infty$ when $\operatorname{gr}\rho$ is injective.
\end{definition}

For the Heisenberg algebra with $\Delta=E$ as above, $\mathfrak v=\mathfrak v_1=E$ and $\operatorname{gr}\rho=\rho$, so $r_*=1$, since $Z\in\ker(\operatorname{gr}\rho)$.

\subsection{The common quotient for a fixed grading}

\begin{proposition}[Common quotient]\label{prop:common-quotient}
Let $W$ be an adapted grading, and let
$\pi_W\colon E\to\mathfrak q_W:=E/J_W$ be the quotient map.
\begin{enumerate}
 \item[(i)] $J_W$ is an ideal of $\mathfrak g$, and
       $[\cdot,\cdot]_{\mathfrak g}$ and $[\cdot,\cdot]_{0,W}$ induce
       the same bracket on $\mathfrak q_W$. With the layers
       $\pi_W(W_j)$, $\mathfrak q_W$ is a Carnot algebra, and $\pi_W$ is
       a Lie algebra homomorphism both on $\mathfrak g$ and on
       $\mathfrak g_0^W$.
 \item[(ii)] For every Carnot algebra $\mathfrak h$ and every Lie algebra
       homomorphism $L\colon\mathfrak g\to\mathfrak h$ with
       $L(W_j)\subseteq\mathfrak h_j$ for every $j$, there is a unique
       graded Lie algebra homomorphism
       $\overline L\colon\mathfrak q_W\to\mathfrak h$ with
       $L=\overline L\circ\pi_W$.
\end{enumerate}
\end{proposition}

\begin{proof}
(i) For $x\in E$ and $z\in J_W$,
\[
 [x,z]_{\mathfrak g}=[x,z]_{0,W}+B_W(x,z)\in J_W\,,
\]
because $J_W$ is an ideal of $\mathfrak g_0^W$ containing every bracket
defect. Hence $J_W$ is an ideal of $\mathfrak g$ and
$\pi_W([u,v]_{\mathfrak g})=\pi_W([u,v]_{0,W})$ for all $u,v\in E$:
the two brackets induce the same bracket on $\mathfrak q_W$ and $\pi_W$
is a homomorphism on both algebras. As $J_W$ is a homogeneous ideal of
the Carnot algebra $\mathfrak g_0^W$ (by Lemma~\ref{lem:osculating}(i)),
the quotient $\mathfrak q_W$ is a Carnot algebra with layers
$\pi_W(W_j)$.

(ii) By Lemma~\ref{lem:gr-adapted}, $L$ is horizontal and
$\operatorname{gr}L=L\circ\sigma_W$, that is,
$L=D\circ\sigma_W^{-1}$ with $D:=\operatorname{gr}L$, a graded Lie
algebra homomorphism by Lemma~\ref{lem:horizontal-gr}. Since $L$
preserves $[\cdot,\cdot]_{\mathfrak g}$, Lemma~\ref{lem:kernel-criterion}
gives $J_W\subseteq\ker L$. Hence $L$ induces a Lie algebra homomorphism
$\overline L\colon\mathfrak q_W\to\mathfrak h$ with
$L=\overline L\circ\pi_W$, which is graded because
$\overline L(\pi_W(W_j))=L(W_j)\subseteq\mathfrak h_j$. It is unique
because $\pi_W$ is surjective.

\end{proof}

Following \cite[Definition~3.18]{LDNGT}, an ideal $I$ of $\mathfrak g$ is
a \emph{Carnot quotient ideal with respect to the grading $W$} if
$\mathfrak g/I$ is a Carnot algebra with layers the images of the
$W_j$ (their definition also requires $\Delta\subseteq W_1+I$,
which holds automatically here since $W_1=\Delta$). By
\cite[Lemma~3.23]{LDNGT}, such an ideal is also an ideal of
$\mathfrak g_0^W$, and the two brackets induce the same bracket on the
quotient. The next corollary shows that
Proposition~\ref{prop:common-quotient} singles out the smallest one,
namely the defect ideal $J_W$.

\begin{corollary}\label{cor:least-compatible}
For every adapted grading $W$, the defect ideal $J_W$ is the smallest
Carnot quotient ideal with respect to $W$.
\end{corollary}

\begin{proof}
By Proposition~\ref{prop:common-quotient}(i), $J_W$ is a Carnot quotient
ideal with respect to $W$. If $I$ is another one, the quotient map
$\pi_I\colon\mathfrak g\to\mathfrak g/I$ is a Lie algebra homomorphism
mapping $W_j$ into the $j$-th layer of $\mathfrak g/I$, so it factors
through $\pi_W$ by Proposition~\ref{prop:common-quotient}(ii). Hence
$J_W=\ker\pi_W\subseteq\ker\pi_I=I$.

\end{proof}

\begin{remark}[The defect at the group level]
\label{rem:group-product-defect}
Fix an adapted grading $W$. For $u,v\in E$, define the product defect
in exponential coordinates by
\[
 \eta_W(u,v)
 := \bigl(-(u\star_{0,W}v)\bigr)\star_{0,W}
    (u\star_{\mathfrak g}v)\,,
\]
so that
\[
 u\star_{\mathfrak g}v
 = (u\star_{0,W}v)\star_{0,W}\eta_W(u,v)\,.
\]

The defect ideal $J_W$ is the smallest homogeneous ideal of
$\mathfrak g_0^W$ containing all values of $\eta_W$. Indeed, by
Proposition~\ref{prop:common-quotient}(i), the two BCH products agree
on $E/J_W$, so $\eta_W(u,v)\in J_W$ for all $u,v\in E$. Conversely,
the BCH formula gives
\[
 \eta_W(\tau u,\tau v)
 = \frac{\tau^2}{2}B_W(u,v)+O(\tau^3)\
 \quad\text{as }\tau\to0\,,
\]
where $\tau u$ and $\tau v$ denote ordinary scalar multiplication.
Thus any homogeneous ideal containing all values of $\eta_W$ also
contains
\[
 B_W(u,v)=\lim_{\tau\to0}\frac{2}{\tau^2}\eta_W(\tau u,\tau v)\,,
\]
since it is a closed linear subspace of $E$. By the definition of
$J_W$, it therefore contains $J_W$.

In particular, every graded homomorphism
$A\colon\mathfrak g_0^W\to\mathfrak h$ satisfies, in exponential
coordinates,
\[
 A(u\star_{\mathfrak g}v)
 = A(u)\star_{\mathfrak h}A(v)
   \star_{\mathfrak h}A\bigl(\eta_W(u,v)\bigr)\,.
\]
Hence $A$ preserves the original group product precisely when
$A(\eta_W(u,v))=0$ for all $u,v\in E$, or equivalently when
$J_W\subseteq\ker A$.
\end{remark}

\subsection{Comparison of the two quotients}\label{subsec:comparison}

Let $W$ be an adapted grading. For $\xi\in\mathfrak v_j$, the vector
$\sigma_W(\xi)\in\mathcal F_jE$ represents $\xi$, so
$\rho(\sigma_W(\xi))\in\mathfrak G_{\leq j}$ by \eqref{eq:rho-strict},
and its component of weight $j$ is $\operatorname{gr}\rho(\xi)$. The
components of lower weight may depend on $W$.

\begin{definition}\label{def:correction}
Let $W$ be an adapted grading. The \emph{correction map} of $W$ is
\[
 t_W:=\rho\circ\sigma_W-\operatorname{gr}\rho
 \colon\mathfrak v\longrightarrow\mathfrak G\ ,
 \quad
 t_W(\mathfrak v_j)\subseteq\mathfrak G_{\leq j-1}\,,
\]
and $M_W$ is the ideal of $\mathfrak G$ generated by the homogeneous
components of the vectors $t_W(\xi)$, $\xi\in\mathfrak v$. We call $W$
\emph{minimal} if $t_W=0$, that is if
$\rho\circ\sigma_W=\operatorname{gr}\rho$.
\end{definition}

Being generated by homogeneous elements, $M_W$ is a homogeneous ideal,
so $\mathfrak G/M_W$ is a Carnot algebra.

\begin{proposition}\label{prop:correction-ideal}
For every adapted grading $W$,
\begin{equation}\label{eq:correction-ideal}
 J_W=\rho^{-1}(M_W)\ \text{ and }
 \
 j_W=(\operatorname{gr}\rho)^{-1}(M_W)\,,
\end{equation}
and $\rho$ induces a graded Lie algebra isomorphism
$\mathfrak q_W\to\mathfrak G/M_W$.
\end{proposition}

\begin{proof}
By Proposition~\ref{prop:common-quotient}(i), $\pi_W$ is a Lie algebra
homomorphism from $\mathfrak g$ into the Carnot algebra $\mathfrak q_W$,
and it is horizontal since $\pi_W(\Delta)=\pi_W(W_1)$ is the first layer.
Proposition~\ref{prop:horizontal-quotient}(ii) gives a graded homomorphism
$\phi_W\colon\mathfrak G\to\mathfrak q_W$ with $\pi_W=\phi_W\circ\rho$.
For $\xi\in\mathfrak v_j$,
\[
 \pi_W(\sigma_W(\xi))
 =\phi_W(\operatorname{gr}\rho(\xi))+\phi_W(t_W(\xi))\,.
\]
The left-hand side and the first term on the right lie in the $j$th layer of $\mathfrak q_W$, while the last term has components of
weight less than $j$ only. Hence $\phi_W(t_W(\xi))=0$. As $\ker\phi_W$
is homogeneous, it contains the homogeneous components of $t_W(\xi)$,
so $M_W\subseteq\ker\phi_W$ and
$\rho^{-1}(M_W)\subseteq\ker(\phi_W\circ\rho)=J_W$.

Conversely, let $\pi_{M_W}\colon\mathfrak G\to\mathfrak G/M_W$ be the quotient
map. Since $t_W(\mathfrak v)\subseteq M_W$, for $\xi\in\mathfrak v_j$ we
get $\pi_{M_W}(\rho(\sigma_W(\xi)))=\pi_{M_W}(\operatorname{gr}\rho(\xi))$, which lies
in the $j$-th layer. Thus the Lie algebra homomorphism
$\pi_{M_W}\circ\rho$ maps $W_j$ into the $j$-th layer of $\mathfrak G/M_W$, and
Proposition~\ref{prop:common-quotient}(ii) gives
$J_W\subseteq\ker(\pi_{M_W}\circ\rho)=\rho^{-1}(M_W)$. This proves the first
equality. For the second, since $t_W(\xi)\in M_W$,
\[
 \xi\in j_W
 \iff\rho(\sigma_W(\xi))=\operatorname{gr}\rho(\xi)+t_W(\xi)\in M_W
 \iff\operatorname{gr}\rho(\xi)\in M_W\,.
\]
Finally, $\pi_{M_W}\circ\rho$ is surjective with kernel $J_W$, so it induces a
Lie algebra isomorphism $\mathfrak q_W\to\mathfrak G/M_W$ and it is graded since $\pi_{M_W}\circ\rho$ maps $W_j$ into the $j$-th
layer.

\end{proof}

\section{Which defect ideals occur}\label{sec:classification}

By Proposition~\ref{prop:correction-ideal}, the defect ideals $J_W$ and
$j_W$ are determined by the homogeneous ideal $M_W\subseteq\mathfrak G$.
This section studies which ideals $M_W$ can occur. We first describe all
adapted gradings relative to a minimal one, then count homogeneous
generators, and finally prove Theorem~\ref{thm:B} and study some explicit examples.

\subsection{Minimal gradings}

An adapted grading $W$ is determined by its map $\sigma_W$, since
$W_j=\sigma_W(\mathfrak v_j)$. Conversely, let
$\sigma\colon\mathfrak v\to E$ be a linear map such that, for every $j$
and every $\xi\in\mathfrak v_j$, the vector $\sigma(\xi)\in\mathcal F_jE$
represents $\xi$. Then $\sigma(\mathfrak v_j)$ is a complement of
$\mathcal F_{j-1}E$ in $\mathcal F_jE$ for every $j$, so
$W_j:=\sigma(\mathfrak v_j)$ is an adapted grading with $\sigma_W=\sigma$.

\begin{lemma}\label{lem:minimal-grading}$\phantom{=}$
\begin{enumerate}
 \item[(i)]  The pair $(\mathfrak g,\Delta)$ admits a minimal adapted grading
       $W^0$, that is, one for which
       $\rho\circ\sigma_{W^0}=\operatorname{gr}\rho$.
 \item[(ii)] Let $W^0$ be a minimal adapted grading. For every adapted grading
       $W$ there is a unique linear map $N\colon\mathfrak v\to\mathfrak v$
       with $N(\mathfrak v_j)\subseteq\mathfrak v_{<j}$ for every $j$ and
       \begin{equation}\label{eq:splitting-parameter}
        \sigma_W=\sigma_{W^0}\circ(\mathrm{id}_{\mathfrak v}+N)\,.
       \end{equation}
       Conversely, every such $N$ arises in this way from a unique
       adapted grading $W$. In both cases
       $t_W=\operatorname{gr}\rho\circ N$.
\end{enumerate}
\end{lemma}

\begin{proof}
(i) Let $W$ be any adapted grading. For $\xi\in\mathfrak v_j$, the
correction $t_W(\xi)$ lies in
$\mathfrak G_{\leq j-1}=\rho(\mathcal F_{j-1}E)$ by
\eqref{eq:rho-strict}. Hence there is a linear map
$a\colon\mathfrak v\to E$ with $a(\mathfrak v_j)\subseteq\mathcal F_{j-1}E$
and $\rho\circ a=t_W$. The map $\sigma:=\sigma_W-a$ still sends each
$\xi\in\mathfrak v_j$ to a representative of $\xi$, so $\sigma=\sigma_{W^0}$
for an adapted grading $W^0$, and
$\rho\circ\sigma_{W^0}=\rho\circ\sigma_W-t_W=\operatorname{gr}\rho$.

(ii) For $\xi\in\mathfrak v_j$, the vectors $\sigma_W(\xi)$ and
$\sigma_{W^0}(\xi)$ both represent $\xi$, so their difference lies in
$\mathcal F_{j-1}E=\sigma_{W^0}(\mathfrak v_{<j})$. Hence
$N:=\sigma_{W^0}^{-1}\circ\sigma_W-\mathrm{id}_{\mathfrak v}$ maps
$\mathfrak v_j$ into $\mathfrak v_{<j}$, and it is the only map
satisfying \eqref{eq:splitting-parameter}. Conversely, if $N$ lowers
weights, then $\sigma_{W^0}(\xi+N\xi)$ represents $\xi$ for every
$\xi\in\mathfrak v_j$, so $\sigma_{W^0}\circ(\mathrm{id}_{\mathfrak v}+N)$
is the map $\sigma_W$ of a unique adapted grading $W$. Finally, since
$\rho\circ\sigma_{W^0}=\operatorname{gr}\rho$,
\[
 t_W=\rho\circ\sigma_{W^0}\circ(\mathrm{id}_{\mathfrak v}+N)
 -\operatorname{gr}\rho=\operatorname{gr}\rho\circ N\,.
\]

\end{proof}

\begin{corollary}\label{cor:least-ideals}
For every adapted grading $W$,
$\ker\rho\subseteq J_W$ and $\ker(\operatorname{gr}\rho)\subseteq j_W$,
and each equality holds if and only if $W$ is minimal. If $W^0$ is
minimal, then the minimal adapted gradings are those given by
\eqref{eq:splitting-parameter} with
$N(\mathfrak v)\subseteq\ker(\operatorname{gr}\rho)$. In particular, the
minimal adapted grading is unique if and only if
$\ker(\operatorname{gr}\rho)\subseteq\mathfrak v_s$.
\end{corollary}

\begin{proof}
The two inclusions follow from \eqref{eq:correction-ideal}. Since $\rho$
and $\operatorname{gr}\rho$ are surjective, each of the equalities
$J_W=\ker\rho$ and $j_W=\ker(\operatorname{gr}\rho)$ is equivalent to
$M_W=\{0\}$, that is, to $t_W=0$. Now fix a minimal adapted grading
$W^0$ and write $\sigma_W=\sigma_{W^0}\circ(\mathrm{id}_{\mathfrak v}+N)$
as in Lemma~\ref{lem:minimal-grading}(ii), so that
$t_W=\operatorname{gr}\rho\circ N$. Hence $W$ is
minimal if and only if $\operatorname{gr}\rho\circ N=0$, that is, if and
only if $N(\mathfrak v)\subseteq\ker(\operatorname{gr}\rho)$.

Finally, suppose that $\ker(\operatorname{gr}\rho)\subseteq\mathfrak v_s$,
and let $W$ be minimal. Then $N$ maps each $\mathfrak v_j$ into
$\mathfrak v_{<j}\cap\ker(\operatorname{gr}\rho)\subseteq
\mathfrak v_{<j}\cap\mathfrak v_s$, which is $\{0\}$ because $j\leq s$.
So $N=0$ and $W=W^0$. Otherwise, since $\ker(\operatorname{gr}\rho)$ is
homogeneous, it contains a nonzero vector $\xi$ of some weight $k<s$.
As $\mathfrak v_s\neq\{0\}$, there is a nonzero linear map
$N\colon\mathfrak v\to\mathfrak v$ that vanishes on $\mathfrak v_{<s}$
and maps $\mathfrak v_s$ into the span of $\xi$. Then
$N(\mathfrak v_j)\subseteq\mathfrak v_{<j}$ for every $j$ and
$N(\mathfrak v)\subseteq\ker(\operatorname{gr}\rho)$, so $N$ defines a
minimal adapted grading different from $W^0$.

\end{proof}

\begin{remark}[Minimal adapted gradings need not be unique]\label{rem:nonunique}
 Let $\mathfrak g$ be the
algebra $N_{6,2,10}$ of \cite{LDT22}, with basis $X_1,\dots,X_6$ and
nonzero brackets
\[
 [X_1,X_2]=X_3\ ,\quad [X_1,X_3]=X_5\ ,\quad
 [X_2,X_3]=X_6\ ,\quad [X_2,X_4]=X_5\,,
\]
and let $\Delta=\operatorname{span}\{X_1,X_2,X_4\}$.

\emph{The filtration.} The brackets of elements of $\Delta$ give
$\mathcal F_2E=\operatorname{span}\{X_1,\dots,X_5\}$ and $\mathcal F_3E=E$, so $s=3$. Writing $\xi_i$
for the class of $X_i$ in $\mathfrak v$, we have
$\mathfrak v_1=\operatorname{span}\{\xi_1,\xi_2,\xi_4\}$,
$\mathfrak v_2=\operatorname{span}\{\xi_3,\xi_5\}$ and $\mathfrak v_3=\operatorname{span}\{\xi_6\}$, with brackets $[\xi_1,\xi_2]=\xi_3$, $[\xi_2,\xi_4]=\xi_5$ and
$[\xi_2,\xi_3]=\xi_6$. Note that $[\xi_1,\xi_3]=0$, since
$X_5\in\mathcal F_2E$.

\emph{The quotient.} Let $L$ be a horizontal homomorphism into a Carnot
algebra $\mathfrak h$, we have
$L(X_5)=[L(X_2),L(X_4)]=[L(X_1),L(X_3)]\in\mathfrak h_2\cap\mathfrak h_3$, hence $L(X_5)=0$. On the other
hand, $E/\operatorname{span}\{X_5\}$ is a Carnot algebra with $X_1,X_2,X_4$ of
weight~$1$, $X_3$ of weight $2$ and $X_6$ of weight $3$: its only
nonzero brackets are $[X_1,X_2]=X_3$ and $[X_2,X_3]=X_6$ (this is the
Engel algebra times $\mathbb R$). Every horizontal $L$ therefore
factors through the quotient map by a homomorphism that is graded,
since it maps the generators $X_1,X_2,X_4$ into $\mathfrak h_1$. By
Corollary~\ref{cor:horizontal-quotient}(ii), $\mathfrak G=E/\operatorname{span} \{X_5\}$
and $\rho$ is the quotient map. Consequently $\operatorname{gr}\rho$
maps $\xi_5$ to $0$ and every other $\xi_i$ to the class of $X_i$, so
$\ker(\operatorname{gr}\rho)=\operatorname{span}\{\xi_5\}$ has weight $2<s$.

\emph{Non-uniqueness.} For $\lambda\in\mathbb R$, let $W^\lambda$ be the
adapted grading with $W^\lambda_1=\Delta$,
$W^\lambda_2=\operatorname{span}\{X_3,X_5\}$ and
$W^\lambda_3=\operatorname{span} \{X_6+\lambda X_5\}$. Its map $\sigma_{W^\lambda}$
sends $\xi_6$ to $X_6+\lambda X_5$ and every other $\xi_i$ to $X_i$.
Since $\rho(X_5)=0$, we get $\rho\circ\sigma_{W^\lambda}=\operatorname{gr}\rho$,
so every $W^\lambda$ is minimal, and $W^\lambda\neq W^\mu$ for
$\lambda\neq\mu$. In the notation of Lemma~\ref{lem:minimal-grading},
$W^\lambda$ corresponds to the map $N$ with $N(\xi_6)=\lambda\xi_5$ and
$N=0$ on the other $\xi_i$, which takes values in
$\ker(\operatorname{gr}\rho)$, as pointed out in Corollary~\ref{cor:least-ideals}. The parameter $\lambda$ is invisible in $\mathfrak G$, but not
in the osculating algebra: $[X_2,X_3]_{0,W^\lambda}=X_6+\lambda X_5$.
\end{remark}

\subsection{Homogeneous generators}

For a homogeneous ideal $M\subseteq\mathfrak G$, write $M_j:=M\cap\mathfrak G_j$
and let $[\mathfrak G,M]$ be the span of the brackets
$[\zeta,m]_{\mathfrak G}$ with $\zeta\in\mathfrak G$ and $m\in M$. By the
Jacobi identity, $[\mathfrak G,M]$ is an ideal contained in $M$, and it
is homogeneous because
$[\mathfrak G,M]_j=\sum_{a+b=j}[\mathfrak G_a,M_b]_{\mathfrak G}$. Thus
\[
 \mu_j(M)=\dim M_j-\dim[\mathfrak G,M]_j\,.
\]
The next lemma shows that $\mu_j(M)$ is the least number of generators
of weight $j$ needed to generate $M$. It is the graded version of Nakayama lemma
\cite[Lemma~1.22]{Rog14} for $M$ as a graded module over the universal
enveloping algebra $U(\mathfrak G)$, acting by the adjoint action: the
elements of positive weight of $U(\mathfrak G)$ map $M$ onto
$[\mathfrak G,M]$. Since that reference works in much greater
generality, we include the direct proof in our setting.

\begin{lemma}\label{lem:minimal-generators}
Let $M\subseteq\mathfrak G$ be a homogeneous ideal, and let
$U_j\subseteq M_j$ be subspaces, $j\geq1$.
\begin{enumerate}
 \item[(i)] If the $U_j$ generate $M$ as an ideal, then
       $M_j=U_j+[\mathfrak G,M]_j$, and so $\mu_j(M)\leq\dim U_j$, for
       every $j$.
 \item[(ii)] If $M_j=U_j\oplus[\mathfrak G,M]_j$ for every $j$, then the $U_j$
       generate $M$ as an ideal and $\dim U_j=\mu_j(M)$.
\end{enumerate}
\end{lemma}

\begin{proof}
(i) The ideal generated by the $U_j$ is spanned by the $U_j$ and by
iterated brackets of their elements with elements of $\mathfrak G$, which belong to $[\mathfrak G,M]$. Hence
$M=\sum_jU_j+[\mathfrak G,M]$, and taking components of weight $j$ gives
the claim.

(ii) Let $I$ be the ideal generated by the $U_j$, then $I\subseteq M$.
Since $M=\bigoplus_jM_j$, it suffices to show that $M_j\subseteq I$ for
every $j$. We argue by induction on $j$, assuming that
$M_b\subseteq I$ for all $b<j$. Let $m\in M_j$, and write $m=u+z$ with
$u\in U_j$ and $z\in[\mathfrak G,M]_j$. Since $\mathfrak G$ and $M$ are
spanned by homogeneous elements,
\[
 [\mathfrak G,M]_j=\sum_{a+b=j}[\mathfrak G_a,M_b]_{\mathfrak G}\,,
 \qquad a,b\geq1\,.
\]
Every $M_b$ occurring here has $b\leq j-1$, so $M_b\subseteq I$ by the
induction hypothesis, and $[\mathfrak G_a,M_b]_{\mathfrak G}\subseteq I$
because $I$ is an ideal. Hence $z\in I$, and $m=u+z\in I$ since
$u\in U_j\subseteq I$. For $j=1$, the equality hodls trivially since the sum is empty. Finally,
$\dim U_j=\dim M_j-\dim[\mathfrak G,M]_j=\mu_j(M)$.

\end{proof}

\subsection{The classification}

\begin{theorem}\label{thm:attainability}
A homogeneous ideal $M\subseteq\mathfrak G$ equals $M_W$ for some
adapted grading $W$ if and only if
\[
 \mu_j(M)\leq\dim\mathfrak v_{>j}\qquad\text{for every }j\geq1\,.
\]
Consequently, as $W$ ranges over all adapted gradings, $J_W$ ranges
exactly over the ideals $\rho^{-1}(M)$, and $j_W$ over the ideals
$(\operatorname{gr}\rho)^{-1}(M)$, with $M$ as above.
\end{theorem}

\begin{proof}
Fix a minimal adapted grading $W^0$, which exists by
Lemma~\ref{lem:minimal-grading}(i). Let $W$ be an adapted grading, and
let $N$ be as in \eqref{eq:splitting-parameter}. Since
$t_W=\operatorname{gr}\rho\circ N$ and $\operatorname{gr}\rho$ is graded,
the component of weight $j$ of $t_W(\xi)$ is
$\operatorname{gr}\rho\bigl(\operatorname{pr}_{\mathfrak v_j}N(\xi)\bigr)$
for every $\xi\in\mathfrak v$. Moreover,
$\operatorname{pr}_{\mathfrak v_j}\circ N$ vanishes on
$\mathfrak v_{\leq j}$, hence, as $\xi$ ranges over $\mathfrak v$, the components
of weight $j$ of $t_W(\xi)$ span $U_j:=\operatorname{gr}\rho\bigl(\operatorname{pr}_{\mathfrak v_j}N(\mathfrak v_{>j})\bigr)$. In particular
$\dim U_j\leq\dim\mathfrak v_{>j}$.
By definition the $U_j$ generate $M_W$, so
Lemma~\ref{lem:minimal-generators}(i) gives
$\mu_j(M_W)\leq\dim\mathfrak v_{>j}$.

Conversely, let $M$ satisfy the inequalities. We construct $N$ so that,
for every $j$, the components of weight $j$ of $t_W$ span a prescribed
generating subspace of $M_j$.

For each $j$, choose a complement $U_j$ of $[\mathfrak G,M]_j$ in $M_j$.
By Lemma~\ref{lem:minimal-generators}(ii), the $U_j$ generate $M$, and
$\dim U_j=\mu_j(M)\leq\dim\mathfrak v_{>j}$. Hence there is a surjective
linear map $T_j\colon\mathfrak v_{>j}\to U_j$. Since
$\operatorname{gr}\rho$ maps $\mathfrak v_j$ onto $\mathfrak G_j$, there is also a linear
map $\ell_j\colon\mathfrak G_j\to\mathfrak v_j$ with
$\operatorname{gr}\rho\circ\ell_j=\mathrm{id}_{\mathfrak G_j}$, i.e. it lifts
$U_j\subseteq\mathfrak G_j$ back to $\mathfrak v_j$. Define
\[
 N(\xi):=\sum_{j\geq1}\ell_j\,T_j\bigl(\operatorname{pr}_{\mathfrak v_{>j}}\xi\bigr)\ \text{ for any }\ \xi\in\mathfrak v\,.
\] If $\xi\in\mathfrak v_k$, then
$\operatorname{pr}_{\mathfrak v_{>j}}\xi$ is $\xi$ for $j<k$ and $0$ for
$j\geq k$, so $N(\xi)=\sum_{j<k}\ell_jT_j(\xi)\in\mathfrak v_{<k}$. Thus
$N$ lowers weights, and by Lemma~\ref{lem:minimal-grading}(ii) it
defines an adapted grading $W$ with $t_W=\operatorname{gr}\rho\circ N$.
For $\xi\in\mathfrak v$, the component of weight $j$ of $N(\xi)$ is
$\ell_jT_j(\operatorname{pr}_{\mathfrak v_{>j}}\xi)$, so the component of
weight $j$ of $t_W(\xi)$ is
$\operatorname{gr}\rho\,\ell_jT_j(\operatorname{pr}_{\mathfrak v_{>j}}\xi)
=T_j(\operatorname{pr}_{\mathfrak v_{>j}}\xi)$. As $\xi$ ranges over
$\mathfrak v$, these span $T_j(\mathfrak v_{>j})=U_j$. Hence $M_W$, the
ideal generated by these components, is the ideal generated by the
$U_j$, that is, $M_W=M$.

The final claim then holds by Proposition~\ref{prop:correction-ideal} as
$J_W=\rho^{-1}(M_W)$ and $j_W=(\operatorname{gr}\rho)^{-1}(M_W)$ for
every adapted grading $W$.

\end{proof}

For $j\geq s$ the condition reads $\mu_j(M)=0$: an attainable ideal needs
no generator of weight $s$, although it may contain elements of weight
$s$.

\subsection{Examples}

\begin{example}\label{ex:heisenberg-all}
Let $\mathfrak g$ and $\Delta$ be as in the Heisenberg example of
the introduction. Every horizontal homomorphism $L$ into a Carnot
algebra is graded, since $L(Z)=[L(X),L(Y)]\in\mathfrak h_2$. Hence
$\mathfrak G=\mathfrak g$ with its stratification, $\rho=\mathrm{id}_E$,
and $\operatorname{gr}\rho=\sigma_{W^0}$ for the stratification given by the unique minimal grading $W^0$. Every adapted grading has
$W_2=\operatorname{span}\{Z+u\}$ for a given $u\in\Delta$, so that
$t_W(\overline Z)=u$ and $t_W$ vanishes on $\mathfrak v_1$. Since
$[u,\Delta]_{\mathfrak g}=\operatorname{span}\{Z\}$ when $u\neq0$,
Proposition~\ref{prop:correction-ideal} gives
\[
 J_W=\{0\}\ \text{ if }u=0\ \text{ and }
 \
 J_W=\operatorname{span}\{u\}\oplus\operatorname{span}\{Z\}\ \text{ if }u\neq0\,.
\]
This agrees with Theorem~\ref{thm:attainability}. A homogeneous ideal of
$\mathfrak g$ with a nonzero component in $\Delta$ must contain $Z$, since
$[u,\Delta]_{\mathfrak g}=\operatorname{span}\{Z\}$ for every
$u\in\Delta$. So the homogeneous ideals of $\mathfrak g$ are exactly
$\{0\}$, $\operatorname{span}\{Z\}$, $\ell\oplus\operatorname{span}\{Z\}$
for a one-dimensional subspace $\ell\subseteq\Delta$, and $\mathfrak g$.
Here $\dim\mathfrak v_{>1}=1$ and $\dim\mathfrak v_{>2}=0$, so the ideals
$\operatorname{span}\{Z\}$ (with $\mu_2=1$) and $\mathfrak g$ (with
$\mu_1=2$) violate the bounds, while every
$\ell\oplus\operatorname{span}\{Z\}$ (with $\mu_1=1$) satisfies them.
Indeed, $\ell\oplus\operatorname{span}\{Z\}=J_W$ for any grading with
$W_2=\operatorname{span}\{Z+u\}$ and $\ell=\operatorname{span}\{u\}$. In particular, the centre $\operatorname{span}\{Z\}$ is never a
defect ideal.
\end{example}

Next, let $\mathfrak g$ be the five-dimensional algebra $N_{5,1}$
\cite{LDT22}, with basis $X_1,\dots,X_5$ and nonzero brackets
\begin{equation}\label{eq:n51-brackets}
 [X_1,X_2]=X_3\ ,\quad [X_1,X_3]=X_4\ ,\quad
[X_2,X_3]=X_5\ ,\quad [X_1,X_4]=X_5\,.
\end{equation}
Let
$\Delta=\operatorname{span}\{X_1,X_2\}$. Then
$\mathcal F_2E=\operatorname{span}\{X_1,X_2,X_3\}$ and
$\mathcal F_3E=E$, so $s=3$, whereas the nilpotency step is $4$. Writing $\xi_i$ for the class of $X_i$ in
$\mathfrak v$, the nonzero brackets of $\mathfrak v$ are
$[\xi_1,\xi_2]=\xi_3$, $[\xi_1,\xi_3]=\xi_4$ and $[\xi_2,\xi_3]=\xi_5$,
hence $\mathfrak v$ is the free nilpotent Lie algebra of rank $2$ and
step $3$. 

Let $\mathfrak e$ be the Engel algebra, with basis
$Y_1,Y_2\mid Y_3\mid Y_4$ listed by weight and nonzero brackets
$[Y_1,Y_2]=Y_3$ and $[Y_1,Y_3]=Y_4$.

\begin{proposition}\label{prop:n51}
For $(N_{5,1},\Delta=\operatorname{span}\{X_1,X_2\})$, we have:
\begin{enumerate}
 \item[(i)] $\mathfrak G=\mathfrak e$, with $\rho(X_i)=Y_i$ for $i\leq4$ and
       $\rho(X_5)=0$. Hence $\ker(\operatorname{gr}\rho)=\operatorname{span}\{\xi_5\}$
       and $r_*=3$.
 \item[(ii)] The grading $W^0$ with $W^0_2=\operatorname{span}\{X_3\}$ and
       $W^0_3=\operatorname{span}\{X_4,X_5\}$ is the unique minimal
       adapted grading.
 \item[(iii)] As $W$ ranges over all adapted gradings, $J_W$ ranges exactly
       over
       \[
        \operatorname{span}\{X_5\}\ ,\quad
        \operatorname{span}\{X_3,X_4,X_5\}\ ,\quad
        \ell\oplus\operatorname{span}\{X_3,X_4,X_5\}\ \
        (\ell\text{ a line in }\Delta)\ \text{ and }\
        E\,.
       \]
\end{enumerate}
\end{proposition}

\begin{proof}
(i) The map $\rho$ in the statement preserves the brackets
\eqref{eq:n51-brackets}: we have $[Y_1,Y_2]=Y_3$ and
$[Y_1,Y_3]=Y_4$, while $[Y_2,Y_3]=[Y_1,Y_4]=0$. So $\rho$ is a horizontal
homomorphism onto $\mathfrak e$ with kernel $\operatorname{span}\{X_5\}$.
Let $L$ be a horizontal homomorphism into a Carnot algebra
$\mathfrak h$. Then $L(X_3)=[L(X_1),L(X_2)]\in\mathfrak h_2$ and
$L(X_4)=[L(X_1),L(X_3)]\in\mathfrak h_3$, so
$L(X_5)=[L(X_2),L(X_3)]=[L(X_1),L(X_4)]\in\mathfrak h_3\cap\mathfrak h_4$, hence $L(X_5)=0$. Therefore
$L=\overline L\circ\rho$ for a unique homomorphism
$\overline L\colon\mathfrak e\to\mathfrak h$, which is graded because
$\overline L(Y_i)=L(X_i)\in\mathfrak h_1$ for $i=1,2$ and $Y_1,Y_2$
generate $\mathfrak e$. By Corollary~\ref{cor:horizontal-quotient}(ii),
$\mathfrak G=\mathfrak e$, with $\rho$ as quotient map.

Finally,
$\operatorname{gr}\rho(\xi_i)=Y_i$ for $i\leq4$ and
$\operatorname{gr}\rho(\xi_5)=0$.
Since $Y_1,\dots,Y_4$ are linearly independent,
$\ker(\operatorname{gr}\rho)=\operatorname{span}\{\xi_5\}$, which has
weight $3$, so $r_*=3$ by Definition~\ref{def:rstar}.

(ii) We have $\rho\circ\sigma_{W^0}(\xi_i)=\rho(X_i)
=\operatorname{gr}\rho(\xi_i)$ for every $i$, so $W^0$ is minimal. It is
unique by Corollary~\ref{cor:least-ideals}, since
$\ker(\operatorname{gr}\rho)\subseteq\mathfrak v_3=\mathfrak v_s$.

(iii) The homogeneous ideals of $\mathfrak e$ are
\[
 \{0\}\ ,\quad \operatorname{span}\{Y_4\}\ ,\quad
 \operatorname{span}\{Y_3,Y_4\}\ ,\quad
 L_1\oplus\operatorname{span}\{Y_3,Y_4\}\ \ (L_1\text{ a line in }\mathfrak e_1)\ \text{ and }\ \mathfrak e\,.
\]
Here $\dim\mathfrak v_{>1}=3$, $\dim\mathfrak v_{>2}=2$ and
$\dim\mathfrak v_{>3}=0$. The ideal $\operatorname{span}\{Y_4\}$ has $\mu_3=1$ and
violates the bounds. For the other nonzero ideals, the only nonzero
$\mu_j$ are $\mu_2=1$ for $\operatorname{span}\{Y_3,Y_4\}$, $\mu_1=1$ for
$L_1\oplus\operatorname{span}\{Y_3,Y_4\}$, and $\mu_1=2$ for
$\mathfrak e$, so they satisfy the bounds. By
Theorem~\ref{thm:attainability}, the ideals $J_W$ are the preimages
under $\rho$ of $\{0\}$, $\operatorname{span}\{Y_3,Y_4\}$,
$L_1\oplus\operatorname{span}\{Y_3,Y_4\}$ and $\mathfrak e$, which are
the ideals in the statement, with $\ell=\Delta\cap\rho^{-1}(L_1)$.
\end{proof}

\begin{remark}\label{rem:n51-table}
The following adapted gradings, all with $W_1=\Delta$, realise the four
types in Proposition~\ref{prop:n51}(iii). Here we are assuming $\ell=\operatorname{span}\{u\}\subset \Delta$.
\[
\begin{array}{c|c|c|c}
 J_W & W_2 & W_3 & \text{nonzero corrections}\\ \hline
 \operatorname{span}\{X_5\} & \operatorname{span}\{X_3\} & \operatorname{span}\{X_4,X_5\} & \text{none}\\
 \operatorname{span}\{X_3,X_4,X_5\}
   & \operatorname{span}\{X_3\} & \operatorname{span}\{X_4+X_3,X_5\} & t_W(\xi_4)=Y_3\\
 \ell\oplus\operatorname{span}\{X_3,X_4,X_5\}
   & \operatorname{span}\{X_3+u\} & \operatorname{span}\{X_4,X_5\} & t_W(\xi_3)=\rho(u)\\
 E & \operatorname{span}\{X_3+X_1\} & \operatorname{span}\{X_4+X_2,X_5\}
   & t_W(\xi_3)=Y_1,\ t_W(\xi_4)=Y_2
\end{array}
\]
In each row, the ideal of $\mathfrak e$ generated by the corrections is
$\{0\}$, $\operatorname{span}\{Y_3,Y_4\}$,
$\operatorname{span}\{\rho(u)\}\oplus\operatorname{span}\{Y_3,Y_4\}$ and $\mathfrak e$,
respectively, so Proposition~\ref{prop:correction-ideal} gives the
displayed $J_W$.
\end{remark}

\section{The intrinsic criterion}\label{sec:intrinsic}

This section proves Theorems~\ref{thm:A} and~\ref{thm:C}. We first
compute the Pansu differentials of horizontal homomorphisms, then prove
the algebraic criterion and its analytic version, and finally revisit
the examples.

\subsection{Pansu differentials of horizontal homomorphisms}

\begin{lemma}\label{lem:hom-differential}
Let $\Phi\colon G\to H$ be a horizontal Lie group homomorphism. Then
$\Phi$ is Lipschitz, it is Pansu differentiable at every $x\in G$ with
respect to every adapted grading $W$, and
\[
 D^W_P\Phi(x)=\operatorname{gr}d\Phi_e\circ\sigma_W^{-1}\,,
 \qquad
 D_P\Phi(x)=\operatorname{gr}d\Phi_e\,.
\]
\end{lemma}

\begin{proof}
Put $L:=d\Phi_e$. In exponential coordinates $\Phi=L$
(Lemma~\ref{lem:linear-homs}), and $L$ maps horizontal curves to
horizontal curves, multiplying their length by at most the norm of
$L|_\Delta\colon\Delta\to\mathfrak h_1$; so $\Phi$ is Lipschitz. Since
$\Phi$ is a homomorphism,
$\Phi(x)^{-1}\star_{\mathfrak h}\Phi(x\star_{\mathfrak g}v)=L(v)$. For
$u=\sum_ju_j$ with $u_j\in W_j$, Lemma~\ref{lem:horizontal-gr}(i) gives
$L(u_j)\in\mathfrak h_{\leq j}$, so
\[
 R^W_{x,t}(u)=\bar\delta_{1/t}L(\delta^W_tu)
 =\sum_j\sum_{k\leq j}t^{\,j-k}\operatorname{pr}_{\mathfrak h_k}L(u_j)
 \longrightarrow\sum_j\operatorname{pr}_{\mathfrak h_j}L(u_j)
 =\operatorname{gr}L\bigl(\sigma_W^{-1}(u)\bigr)
\]
uniformly on compact sets as $t\downarrow0$, hence also for $d_H$. The
limit $\operatorname{gr}L\circ\sigma_W^{-1}$ is a graded Lie algebra
homomorphism $\mathfrak g_0^W\to\mathfrak h$, by
Lemmata~\ref{lem:horizontal-gr}(ii) and~\ref{lem:osculating}(i).
\end{proof}

\subsection{The algebraic criterion}

\begin{theorem}\label{thm:intrinsic}
Let $D\colon\mathfrak v\to\mathfrak h$ be a graded Lie algebra
homomorphism. The following are equivalent:
\begin{enumerate}
 \item[(i)] $D=\operatorname{gr}d\Phi_e$ for a horizontal Lie group
       homomorphism $\Phi\colon G\to H$;
 \item[(ii)] $D\circ\sigma_W^{-1}$ preserves $[\cdot,\cdot]_{\mathfrak g}$
       for some adapted grading $W$;
 \item[(iii)] $D\circ\sigma_W^{-1}$ preserves $[\cdot,\cdot]_{\mathfrak g}$
       for every minimal adapted grading $W$;
 \item[(iv)] $D(\ker\operatorname{gr}\rho)=\{0\}$;
 \item[(v)] $D=\overline D\circ\operatorname{gr}\rho$ for a graded Lie
       algebra homomorphism $\overline D\colon\mathfrak G\to\mathfrak h$.
\end{enumerate}
If this is the case, $\Phi$ and $\overline D$ are unique,
$d\Phi_e=\overline D\circ\rho$, and $D\circ\sigma_W^{-1}=d\Phi_e$ for
every adapted grading $W$ as in \textup{(ii)}. Moreover,
if $W^0$ is a minimal adapted grading and
$\sigma_W=\sigma_{W^0}\circ(\mathrm{id}_{\mathfrak v}+N)$ as in
Lemma~\ref{lem:minimal-grading}, then
\begin{equation}\label{eq:good-gradings}
 D\circ\sigma_W^{-1}\text{ preserves }[\cdot,\cdot]_{\mathfrak g}
 \ \Longleftrightarrow\ D\circ N=0\,.
\end{equation}
\end{theorem}

\begin{proof}
Put $A_W:=D\circ\sigma_W^{-1}$. By Lemma~\ref{lem:kernel-criterion},
$A_W$ preserves $[\cdot,\cdot]_{\mathfrak g}$ if and only if
$D(j_W)=\{0\}$.

(i)$\Rightarrow$(v). The homomorphism $L:=d\Phi_e$ is horizontal, so
Proposition~\ref{prop:horizontal-quotient}(ii) gives
$D=\operatorname{gr}L=\overline L\circ\operatorname{gr}\rho$.

(v)$\Rightarrow$(iv) is trvial.

(iv)$\Rightarrow$(iii). If $W$ is minimal, then
$j_W=\ker(\operatorname{gr}\rho)$ by Corollary~\ref{cor:least-ideals},
so $D(j_W)=\{0\}$.

(iii)$\Rightarrow$(ii). Minimal adapted gradings exist by
Lemma~\ref{lem:minimal-grading}.

(ii)$\Rightarrow$(i). By Lemma~\ref{lem:kernel-criterion}, $A_W$ is a
horizontal Lie algebra homomorphism with $\operatorname{gr}A_W=D$. Since
$G$ is simply connected, $A_W=d\Phi_e$ for a Lie group homomorphism
$\Phi\colon G\to H$, which is horizontal.

If these hold, $\Phi$ is unique by Lemma~\ref{lem:horizontal-gr}(iii),
since $\operatorname{gr}d\Phi_e=D$ determines $d\Phi_e$, and for the same
reason $A_W=d\Phi_e$ for every $W$ as in (ii). The map $\overline D$ is
unique because $\operatorname{gr}\rho$ is surjective, and so
$\overline D=\overline L$ for $L=d\Phi_e$, that is,
$d\Phi_e=\overline D\circ\rho$.

Finally, by Proposition~\ref{prop:correction-ideal} and the
surjectivity of $\operatorname{gr}\rho$,
\[
 D(j_W)=\overline D\bigl(\operatorname{gr}\rho\bigl((\operatorname{gr}\rho)^{-1}(M_W)\bigr)\bigr)=\overline D(M_W)\,.
\]
Since $\ker\overline D$ is a homogeneous ideal and $M_W$ is generated by
the homogeneous components of the vectors $t_W(\xi)$, this vanishes if
and only if $\overline D\circ t_W=0$. By
Lemma~\ref{lem:minimal-grading},
$\overline D\circ t_W=\overline D\circ\operatorname{gr}\rho\circ N=D\circ N$.
\end{proof}

Let $D$ satisfy the conditions of Theorem~\ref{thm:intrinsic}, and fix
a minimal grading $W^0$. By Lemma~\ref{lem:minimal-grading}(ii), every
adapted grading $W$ is given by
$\sigma_W=\sigma_{W^0}\circ(\mathrm{id}_{\mathfrak v}+N)$ for a unique
linear map $N$ with $N(\mathfrak v_j)\subseteq\mathfrak v_{<j}$ for
every $j$, and by \eqref{eq:good-gradings}, $D\circ\sigma_W^{-1}$
preserves $[\cdot,\cdot]_{\mathfrak g}$ if and only if $N$ takes values
in $\ker D$. Consequently, the larger $\ker D$ is, the more adapted
gradings have this property; if $D$ is injective, $W^0$ is the only one.
Furthermore, the minimal gradings are exactly the adapted gradings $W$
for which $D\circ\sigma_W^{-1}$ preserves $[\cdot,\cdot]_{\mathfrak g}$
for every such $D$, into every Carnot algebra. Indeed, if $W$ is
minimal, then $N$ takes values in $\ker(\operatorname{gr}\rho)$
(Corollary~\ref{cor:least-ideals}), which is contained in $\ker D$.
Conversely, for the target $\mathfrak G$ and $D=\operatorname{gr}\rho$,
the condition reads $N(\mathfrak v)\subseteq\ker(\operatorname{gr}\rho)$,
that is, $W$ is minimal. So a minimal grading can be fixed once and for
all, independently of the target and of $D$, and this is what makes it
the natural choice in Theorem~\ref{thm:intrinsic-analytic}.

\subsection{Pansu differentials}

\begin{theorem}\label{thm:intrinsic-analytic}
Let $f\colon\Omega\to H$ be Lipschitz and $x\in\Omega_f$. The following
are equivalent:
\begin{enumerate}
 \item[(i)] for some adapted grading $W$, $D^W_Pf(x)$ preserves
       $[\cdot,\cdot]_{\mathfrak g}$, that is, $D^W_Pf(x)\colon G\to H$ is a group homomorphism;
 \item[(ii)] $D^W_Pf(x)\colon G\to H$ is a group homomorphism for every minimal adapted grading $W$;
 \item[(iii)] $D_Pf(x)(\ker\operatorname{gr}\rho)=\{0\}$;
 \item[(iv)] there is a horizontal Lie group homomorphism $\Phi\colon G\to H$
       with $D_P\Phi(y)=D_Pf(x)$ for every $y\in G$.
\end{enumerate}
In that case $\Phi$ is unique, and $D^W_Pf(x)=d\Phi_e$ for every
adapted grading $W$ as in \textup{(i)}.
\end{theorem}

\begin{proof}
Set $D:=D_Pf(x)$, we get $D^W_Pf(x)=D\circ\sigma_W^{-1}$ for every $W$
(see Definition~\ref{def:canonical-differential}). By
Lemma~\ref{lem:linear-homs}, a linear map $E\to\mathfrak h$ is a group
homomorphism $G\to H$ if and only if it preserves
$[\cdot,\cdot]_{\mathfrak g}$. By Lemma~\ref{lem:hom-differential}, (iv)
says that $D=\operatorname{gr}d\Phi_e$. So the statements from (i) to (iv) are conditions
(ii), (iii), (iv) and (i) of Theorem~\ref{thm:intrinsic} respectively.

\end{proof}

\begin{corollary}\label{cor:visibility}
Let $m$ be the nilpotency step of $H$ and $\Omega\subseteq G$ open. If $m<r_*$, then for every Lipschitz map
$f\colon\Omega\to H$ the conditions of
Theorem~\ref{thm:intrinsic-analytic} hold at every $x\in\Omega_f$. In
particular, every minimal adapted grading $W$ satisfies
$D^W_Pf(x)=d\Phi_e$ for a horizontal homomorphism $\Phi\colon G\to H$.
Here $\Phi$ depends on $f$ and $x$, but $W$ does not.
\end{corollary}

\begin{proof}
Let $x\in\Omega_f$ and put $D:=D_Pf(x)\colon\mathfrak v\to\mathfrak h$. If $r_*=\infty$, then
$\ker(\operatorname{gr}\rho)=\{0\}$ and there is nothing to prove.
Otherwise, by Definition~\ref{def:rstar} we have
\[
 \ker(\operatorname{gr}\rho)\subseteq\bigoplus_{j\geq r_*}\mathfrak v_j\,.
\]
For $j\geq r_*>m$, the graded map $D$ sends $\mathfrak v_j$ into
$\mathfrak h_j=\{0\}$. Therefore $D(\ker(\operatorname{gr}\rho))=\{0\}$, which yields the claim by condition (iii) of Theorem~\ref{thm:intrinsic-analytic}.

\end{proof}
\begin{remark}\label{rem:graded-targets}
Theorem~\ref{thm:intrinsic} holds verbatim when $\mathfrak h$ is only
graded, since its proof uses no other property of $\mathfrak h$. If $H$
is a graded group with a homogeneous distance, then
Lemma~\ref{lem:change-grading} and Theorem~\ref{thm:intrinsic-analytic}
hold at every point where $f$ is Pansu differentiable, and
Corollary~\ref{cor:visibility} holds with $m$ the step of the Carnot
subalgebra generated by $\mathfrak h_1$, which contains the image of
every Pansu differential. What is missing is the almost everywhere
existence of Pansu differentials of Lipschitz maps from $G$ into $H$,
which combines the non-Carnot domain of Theorem~\ref{thm:existence} with the graded setting of \cite[Theorem~2.11]{Mag10}.
\end{remark}
\subsection{Examples}

\begin{example}\label{ex:heisenberg-revisited}
For the Heisenberg algebra in Example~\ref{ex:heisenberg-all}, we have that $\ker(\operatorname{gr}\rho)=\{0\}$, that is $r_*=\infty$ and
Corollary~\ref{cor:visibility} applies to every Carnot group $H$: for
the stratification $W^0$, $D^{W^0}_Pf(x)$ is a group homomorphism at
every $x\in\Omega_f$. For the identity map of Example~\ref{ex:heisenberg-differentials},
$D=D_P\mathrm{id}(x)=\sigma_{W^0}$ is injective, so by
\eqref{eq:good-gradings} the stratification is the only adapted grading
for which $D^W_P\mathrm{id}(x)$ is a group homomorphism. The failure for
the tilted grading concerns only that particular choice.
\end{example}

The remaining examples concern the bound $m<r_*$ used in
Corollary~\ref{cor:visibility}. Proposition~\ref{prop:n51-engel} shows
that it is sharp in general, while Proposition~\ref{prop:rigidity} shows
that, for a particular pair $(\mathfrak g,\Delta)$, the conclusion of the
corollary may persist when $m=r_*$.

\begin{proposition}\label{prop:n51-engel}
Let $\mathfrak g=N_{5,1}$ and $\Delta=\operatorname{span}\{X_1,X_2\}$ be
as in Proposition~\ref{prop:n51}, so that $r_*=3$, and let
$H=(\mathfrak e,\star_{\mathfrak e})$ be the Engel group, of step $m=3$.
There is a polynomial contact map $F\colon G\to H$ such that, for every
ball $\Omega\subseteq G$, the restriction of $F$ to $\Omega$ is
Lipschitz and
\[
 D_PF(x)(\xi_5)=-15\,Y_4\neq0\qquad\text{for every }x\in\Omega_F\,.
\]
Hence at no point of $\Omega_F$ do the conditions of
Theorem~\ref{thm:intrinsic-analytic} hold: for no adapted grading $W$ is
$D^W_PF(x)$ a group homomorphism, and $D_PF(x)$ is not the Pansu
differential of any horizontal homomorphism.
\end{proposition}

\begin{proof}
In exponential coordinates $x=(x_1,\dots,x_5)$ of $G$, the
left-invariant vector fields of $X_1$ and $X_2$ are
\[
\begin{aligned}
 X_1&=\partial_1-\tfrac{x_2}{2}\partial_3
 -\Bigl(\tfrac{x_3}{2}+\tfrac{x_1x_2}{12}\Bigr)\partial_4
 -\Bigl(\tfrac{x_4}{2}+\tfrac{x_1x_3}{12}+\tfrac{x_2^2}{12}\Bigr)\partial_5\ \text{ and }\ 
 X_2=\partial_2+\tfrac{x_1}{2}\partial_3+\tfrac{x_1^2}{12}\partial_4
 +\Bigl(\tfrac{x_1x_2}{12}-\tfrac{x_3}{2}\Bigr)\partial_5\,,
\end{aligned}
\]
and in exponential coordinates $y=(y_1,\dots,y_4)$ of $H$, with respect
to the basis $Y_1,\dots,Y_4$, those of $Y_1$ and $Y_2$ are
\[
 Y_1=\partial_{y_1}-\tfrac{y_2}{2}\partial_{y_3}
 -\Bigl(\tfrac{y_3}{2}+\tfrac{y_1y_2}{12}\Bigr)\partial_{y_4}\ \text{ and }
 \ 
 Y_2=\partial_{y_2}+\tfrac{y_1}{2}\partial_{y_3}+\tfrac{y_1^2}{12}\partial_{y_4}\,.
\]
A smooth map $F=(F_1,\dots,F_4)\colon G\to H$ is contact, that is,
$dF$ maps $\Delta$ into the horizontal bundle of $H$, if and only if for $i=1,2$
\[
\begin{aligned}
 X_iF_3&=\tfrac12\bigl(F_1\,X_iF_2-F_2\,X_iF_1\bigr)\,,\\
 X_iF_4&=\tfrac{1}{12}F_1^2\,X_iF_2
 -\Bigl(\tfrac12F_3+\tfrac{1}{12}F_1F_2\Bigr)X_iF_1\,.
\end{aligned}
\]
A direct computation shows that these identities hold for the
polynomial map
\[
\begin{aligned}
 F_1&=x_1+x_2+\tfrac12x_1^2\,,\, F_2=15x_1-25\,,\, F_3=-\tfrac54x_1^3+\tfrac{25}{4}x_1^2+\tfrac{25}{2}x_1
 +\tfrac{25}{2}x_2-15x_3-25\,,\\
 F_4&=\tfrac{1}{16}x_1^5-\tfrac{25}{48}x_1^4-\tfrac{25}{12}x_1^3
 -\tfrac{25}{12}x_1^2x_2+\tfrac54x_1^2x_3+\tfrac{25}{6}x_1^2
 -\tfrac{25}{6}x_1x_2-\tfrac{15}{2}x_1x_4\\
 &\quad+\tfrac{25}{2}x_1-\tfrac{25}{12}x_2^2+\tfrac{25}{2}x_2
 -15x_4-15x_5\,.
\end{aligned}
\]
Since $F$ is smooth and contact, it maps horizontal curves to horizontal
curves, with speed bounded on compact sets. Hence it is Lipschitz on
every ball $\Omega$, and so Pansu differentiable almost everywhere on
$\Omega$ by Theorem~\ref{thm:existence}. When the domain is a
stratified group, locally Lipschitz maps are characterised by the
contact equations of \cite[Theorem~1.1]{Mag10} (see also \cite{Mag13}).

Let $x\in\Omega_F$. For $i=1,2$,  $\delta^W_t(X_i)=tX_i$ and left translations
act on the coordinates $y_i$ by addition, so that the
components of weight $1$ of $R^W_{x,t}(X_i)$ are
$\tfrac1t\bigl(F_k(x\star_{\mathfrak g}tX_i)-F_k(x)\bigr)$, $k=1,2$,
which tend to $X_iF_k(x)$. As $X_1F_1=1+x_1$, $X_2F_1=1$, $X_1F_2=15$
and $X_2F_2=0$, this gives
\[
 D_PF(x)(\xi_1)=(1+x_1)\,Y_1+15\,Y_2\ \text{ and }\
 D_PF(x)(\xi_2)=Y_1\,.
\]
Since $D_PF(x)$ is a Lie algebra homomorphism of $\mathfrak v$ and
$\xi_5=[\xi_2,[\xi_1,\xi_2]]$,
\[
 D_PF(x)(\xi_5)
 =\bigl[Y_1,[(1+x_1)Y_1+15\,Y_2,\,Y_1]\bigr]
 =[Y_1,-15\,Y_3]=-15\,Y_4\,.
\]
The claim then follows from Theorem~\ref{thm:intrinsic-analytic},
since $\xi_5\in\ker(\operatorname{gr}\rho)$.

\end{proof}

\begin{example}\label{ex:rigidity-algebra}
Let $\mathfrak g$ have basis $X_1,\dots,X_6$ and nonzero brackets
\[
 [X_1,X_2]=X_5\ ,\quad [X_1,X_5]=X_6\ ,\quad [X_3,X_4]=X_6\,,
\]
 and let $\Delta=\operatorname{span}\{X_1,\dots,X_4\}$.
Then $\mathcal F_2E=E$, so $s=2$, and the nonzero brackets of
$\mathfrak v$ are $[\xi_1,\xi_2]=\xi_5$ and $[\xi_3,\xi_4]=\xi_6$; in
particular $[\xi_1,\xi_5]=0$, since $X_6\in\mathcal F_2E$. A horizontal
homomorphism $L$ into a Carnot algebra $\mathfrak h$ satisfies
$L(X_6)=[L(X_3),L(X_4)]\in\mathfrak h_2$ and
$L(X_6)=[L(X_1),[L(X_1),L(X_2)]]\in\mathfrak h_3$, so $L(X_6)=0$.
Moreover, $E/\operatorname{span}\{X_6\}$ is a Carnot algebra, with
$X_1,\dots,X_4$ of weight $1$ and $X_5$ of weight $2$: it is the product
of the Heisenberg algebra spanned by $X_1,X_2,X_5$ with the abelian
algebra $\operatorname{span}\{X_3,X_4\}$. As in
Proposition~\ref{prop:n51}(i), Corollary~\ref{cor:horizontal-quotient}(ii)
gives $\mathfrak G=E/\operatorname{span}\{X_6\}$. Hence
$\ker(\operatorname{gr}\rho)=\operatorname{span}\{\xi_6\}$ and $r_*=2$,
and the grading $W^0$ with $W^0_2=\operatorname{span}\{X_5,X_6\}$ is the
unique minimal adapted grading (Corollary~\ref{cor:least-ideals}).

Let $H=H^1$ be the Heisenberg group, with $\mathfrak h$ spanned by
$A,B$ and $[A,B]$. Then $m=2=r_*$, so Corollary~\ref{cor:visibility}
does not apply, and the criterion of Theorem~\ref{thm:intrinsic} is a
genuine condition: the graded homomorphism
$D\colon\mathfrak v\to\mathfrak h$ that maps $\xi_3,\xi_4,\xi_6$ to
$A,B,[A,B]$ and $\xi_1,\xi_2,\xi_5$ to $0$ satisfies $D(\xi_6)\neq0$.
Nevertheless, the next proposition shows that the criterion holds at
every point for maps whose Pansu differential exists everywhere and is
continuous.
\end{example}

\begin{proposition}\label{prop:rigidity}
In Example~\ref{ex:rigidity-algebra}, let $f\colon\Omega\to H$ be
Pansu differentiable with respect to $W^0$ at every point of $\Omega$,
with $x\mapsto D^{W^0}_Pf(x)$ continuous. Then, at every $x\in\Omega$,
$D^{W^0}_Pf(x)$ preserves $[\cdot,\cdot]_{\mathfrak g}$, that is,
$D^{W^0}_Pf(x)=d\Phi_e$ for a unique horizontal homomorphism
$\Phi\colon G\to H$.
\end{proposition}

\begin{proof}
Write $D_x:=D^{W^0}_Pf(x)$. Since $W^0$ is minimal,
$J_{W^0}=\ker\rho=\operatorname{span}\{X_6\}$, so by
Lemma~\ref{lem:kernel-criterion} $D_x$ preserves
$[\cdot,\cdot]_{\mathfrak g}$ if and only if $D_x(X_6)=0$. We argue by
contradiction: suppose that the set
$O:=\{x\in\Omega:D_x(X_6)\neq0\}$ is nonempty. We know it is open, because
$x\mapsto D_x$ is continuous. We recall that in this particular example we are taking $H=H^1$ to be the three-dimensional Heisenberg group. The idea is that on $O$ the map $f$ has zero derivative along $X_1$ and
$X_2$, hence is constant along their flows, and therefore also along
$X_6=[X_1,[X_1,X_2]]$, which contradicts the assumption that $D_x(X_6)\neq0$.

\emph{Step 1: $D_x(X_1)=D_x(X_2)=0$ for $x\in O$.} The vectors
$D_x(X_3)$ and $D_x(X_4)$ span $\mathfrak h_1$, since
$[D_x(X_3),D_x(X_4)]=D_x(X_6)\neq0$. For $i\in\{1,2\}$ and
$j\in\{3,4\}$ we have $[X_i,X_j]_{0,W^0}=0$, so
$[D_x(X_i),D_x(X_j)]=0$. Since the bracket
$\mathfrak h_1\times\mathfrak h_1\to\mathfrak h_2$ is nondegenerate, this implies that
$D_x(X_1)=D_x(X_2)=0$.

\emph{Step 2: $f$ is constant along $X_1$- and $X_2$-segments in $O$.}
Let $y\in O$ and $i\in\{1,2\}$. Since $X_i\in W^0_1$, we have
$\delta^{W^0}_t(\pm X_i)=\pm tX_i$. By left invariance and homogeneity
of $d_H$, Definition~\ref{def:pansu-W} and Step~1 give
\[
 \tfrac1t\,d_H\bigl(f(y),f(y\star_{\mathfrak g}(\pm tX_i))\bigr)
 =d_H\bigl(0,R^{W^0}_{y,t}(\pm X_i)\bigr)
 \xrightarrow[t\to 0]{} d_H\bigl(0,D_y(\pm X_i)\bigr)=0\,.
\]
 Now let $I$ be an interval with
$c(s):=y\star_{\mathfrak g}sX_i\in O$ for all $s\in I$. Since $sX_i$ and
$hX_i$ commute, $c(s+h)=c(s)\star_{\mathfrak g}hX_i$, so the estimate
above at the point $c(s)$ gives
$d_H\bigl(f(c(s)),f(c(s+h))\bigr)=o(|h|)$ as $h\to0$, for every
$s\in I$. Fix $s_0\in I$ and put
$\varphi(s):=d_H\bigl(f(c(s_0)),f(c(s))\bigr)$. By the triangle
inequality, $|\varphi(s+h)-\varphi(s)|\leq d_H\bigl(f(c(s)),f(c(s+h))\bigr)$,
so $\varphi$ is differentiable on $I$ with $\varphi'=0$. Hence $\varphi$
is constant, equal to $\varphi(s_0)=0$, and $f\circ c$ is constant on
$I$.

\emph{Step 3: $f(x\star_{\mathfrak g}a^3X_6)=f(x)$ for $x\in O$ and
$|a|$ small.} For $a\in\mathbb R$ put $A:=aX_1$ and $B:=aX_2$, and write
$[g,h]:=g\star_{\mathfrak g}h\star_{\mathfrak g}g^{-1}\star_{\mathfrak g}h^{-1}$
for the group commutator. The BCH formula gives
\[
 [A,B]=a^2X_5+\tfrac12a^3X_6\ \text{ and }\ 
 [A,[A,B]]=a^3X_6\,.
\]
Expanding the commutators, $a^3X_6=g_1\star_{\mathfrak g}\dots\star_{\mathfrak g}g_{10}$
with each $g_k\in\{\pm aX_1,\pm aX_2\}$. Put $x_0:=x$ and
$x_k:=x_{k-1}\star_{\mathfrak g}g_k$, so that
$x_{10}=x\star_{\mathfrak g}a^3X_6$, and note that $x_{k-1}$ and $x_k$
are joined by the segment
$\{x_{k-1}\star_{\mathfrak g}\tau g_k:\tau\in[0,1]\}$ along $X_1$ or
$X_2$. As $a\to0$, these ten segments shrink to the point $x$, so they
lie in the open set $O$ for $|a|$ small. By Step~2,
$f(x_k)=f(x_{k-1})$ for every $k$, and therefore
$f(x\star_{\mathfrak g}a^3X_6)=f(x)$.

\emph{Step 4: conclusion.} Let $x\in O$. Since $X_6\in W^0_2$, we have
$\delta^{W^0}_t(X_6)=t^2X_6$. For small $t>0$, Step~3 with $a=t^{2/3}$
gives $f(x\star_{\mathfrak g}t^2X_6)=f(x)$, so
\[
 R^{W^0}_{x,t}(X_6)
 =\bar\delta_{1/t}\bigl(f(x)^{-1}\star_{\mathfrak h}f(x\star_{\mathfrak g}t^2X_6)\bigr)=0\,.
\]
Letting $t\to0$ gives $D_x(X_6)=0$, which contradicts $x\in O$.

Hence $O=\emptyset$, that is, $D_x(X_6)=0$ for every $x\in\Omega$, and
the last assertion follows from Theorem~\ref{thm:intrinsic}.

\end{proof}

We note that the proof uses that $O$ is open. For a Lipschitz map, the set where
$D_Pf(x)(\xi_6)\neq0$ is only measurable, and we do not know whether the
conclusion of Proposition~\ref{prop:rigidity} holds almost everywhere.


\end{document}